

\ifx\shlhetal\undefinedcontrolsequence\let\shlhetal\relax\fi

\input amstex
\expandafter\ifx\csname mathdefs.tex\endcsname\relax
  \expandafter\gdef\csname mathdefs.tex\endcsname{}
\else \message{Hey!  Apparently you were trying to
  \string\input{mathdefs.tex} twice.   This does not make sense.} 
\errmessage{Please edit your file (probably \jobname.tex) and remove
any duplicate ``\string\input'' lines}\endinput\fi




\catcode`\X=12\catcode`\@=11

\def\n@wcount{\alloc@0\count\countdef\insc@unt}
\def\n@wwrite{\alloc@7\write\chardef\sixt@@n}
\def\n@wread{\alloc@6\read\chardef\sixt@@n}
\def\r@s@t{\relax}\def\v@idline{\par}\def\@mputate#1/{#1}
\def\l@c@l#1X{\firstpart.#1}\def\gl@b@l#1X{#1}\def\t@d@l#1X{{}}

\def\crossrefs#1{\ifx\all#1\let\tr@ce=\all\else\def\tr@ce{#1,}\fi
   \n@wwrite\cit@tionsout\openout\cit@tionsout=\jobname.cit 
   \write\cit@tionsout{\tr@ce}\expandafter\setfl@gs\tr@ce,}
\def\setfl@gs#1,{\def\@{#1}\ifx\@\empty\let\next=\relax
   \else\let\next=\setfl@gs\expandafter\xdef
   \csname#1tr@cetrue\endcsname{}\fi\next}
\def\m@ketag#1#2{\expandafter\n@wcount\csname#2tagno\endcsname
     \csname#2tagno\endcsname=0\let\tail=\all\xdef\all{\tail#2,}
   \ifx#1\l@c@l\let\tail=\r@s@t\xdef\r@s@t{\csname#2tagno\endcsname=0\tail}\fi
   \expandafter\gdef\csname#2cite\endcsname##1{\expandafter
     \ifx\csname#2tag##1\endcsname\relax?\else\csname#2tag##1\endcsname\fi
     \expandafter\ifx\csname#2tr@cetrue\endcsname\relax\else
     \write\cit@tionsout{#2tag ##1 cited on page \folio.}\fi}
   \expandafter\gdef\csname#2page\endcsname##1{\expandafter
     \ifx\csname#2page##1\endcsname\relax?\else\csname#2page##1\endcsname\fi
     \expandafter\ifx\csname#2tr@cetrue\endcsname\relax\else
     \write\cit@tionsout{#2tag ##1 cited on page \folio.}\fi}
   \expandafter\gdef\csname#2tag\endcsname##1{\expandafter
      \ifx\csname#2check##1\endcsname\relax
      \expandafter\xdef\csname#2check##1\endcsname{}%
      \else\immediate\write16{Warning: #2tag ##1 used more than once.}\fi
      \multit@g{#1}{#2}##1/X%
      \write\t@gsout{#2tag ##1 assigned number \csname#2tag##1\endcsname\space
      on page \number\count0.}%
   \csname#2tag##1\endcsname}}

\def\multit@g#1#2#3/#4X{\def\t@mp{#4}\ifx\t@mp\empty%
      \global\advance\csname#2tagno\endcsname by 1 
      \expandafter\xdef\csname#2tag#3\endcsname
      {#1\number\csname#2tagno\endcsnameX}%
   \else\expandafter\ifx\csname#2last#3\endcsname\relax
      \expandafter\n@wcount\csname#2last#3\endcsname
      \global\advance\csname#2tagno\endcsname by 1 
      \expandafter\xdef\csname#2tag#3\endcsname
      {#1\number\csname#2tagno\endcsnameX}
      \write\t@gsout{#2tag #3 assigned number \csname#2tag#3\endcsname\space
      on page \number\count0.}\fi
   \global\advance\csname#2last#3\endcsname by 1
   \def\t@mp{\expandafter\xdef\csname#2tag#3/}%
   \expandafter\t@mp\@mputate#4\endcsname
   {\csname#2tag#3\endcsname\lastpart{\csname#2last#3\endcsname}}\fi}
\def\t@gs#1{\def\all{}\m@ketag#1e\m@ketag#1s\m@ketag\t@d@l p
\let\realscite\scite
\let\realstag\stag
   \m@ketag\gl@b@l r \n@wread\t@gsin
   \openin\t@gsin=\jobname.tgs \re@der \closein\t@gsin
   \n@wwrite\t@gsout\openout\t@gsout=\jobname.tgs }
\outer\def\localtags{\t@gs\l@c@l}
\outer\def\globaltags{\t@gs\gl@b@l}
\outer\def\newlocaltag#1{\m@ketag\l@c@l{#1}}
\outer\def\newglobaltag#1{\m@ketag\gl@b@l{#1}}

\newif\ifpr@ 
\def\m@kecs #1tag #2 assigned number #3 on page #4.%
   {\expandafter\gdef\csname#1tag#2\endcsname{#3}
   \expandafter\gdef\csname#1page#2\endcsname{#4}
   \ifpr@\expandafter\xdef\csname#1check#2\endcsname{}\fi}
\def\re@der{\ifeof\t@gsin\let\next=\relax\else
   \read\t@gsin to\t@gline\ifx\t@gline\v@idline\else
   \expandafter\m@kecs \t@gline\fi\let \next=\re@der\fi\next}
\def\pretags#1{\pr@true\pret@gs#1,,}
\def\pret@gs#1,{\def\@{#1}\ifx\@\empty\let\n@xtfile=\relax
   \else\let\n@xtfile=\pret@gs \openin\t@gsin=#1.tgs \message{#1} \re@der 
   \closein\t@gsin\fi \n@xtfile}

\newcount\sectno\sectno=0\newcount\subsectno\subsectno=0
\newif\ifultr@local \def\ultralocal{\ultr@localtrue}
\def\firstpart{\number\sectno}
\def\lastpart#1{\ifcase#1 \or a\or b\or c\or d\or e\or f\or g\or h\or 
   i\or k\or l\or m\or n\or o\or p\or q\or r\or s\or t\or u\or v\or w\or 
   x\or y\or z \fi}

\def\resetall{\global\advance\sectno by 1\subsectno=0
   \gdef\firstpart{\number\sectno}\r@s@t}
\def\resetsub{\global\advance\subsectno by 1
   \gdef\firstpart{\number\sectno.\number\subsectno}\r@s@t}
\def\newsection#1\par{\resetall\vskip0pt plus.3\vsize\penalty-250
   \vskip0pt plus-.3\vsize\bigskip\bigskip
   \message{#1}\leftline{\bf#1}\nobreak\bigskip}
\def\subsection#1\par{\ifultr@local\resetsub\fi
   \vskip0pt plus.2\vsize\penalty-250\vskip0pt plus-.2\vsize
   \bigskip\smallskip\message{#1}\leftline{\bf#1}\nobreak\medskip}


\newdimen\marginshift

\newdimen\margindelta
\newdimen\marginmax
\newdimen\marginmin

\def\margininit{       
\marginmax=3 true cm                  
				      
\margindelta=0.1 true cm              
\marginmin=0.1true cm                 
\marginshift=\marginmin
}    

\def\t@gsjj#1,{\def\@{#1}\ifx\@\empty\let\next=\relax\else\let\next=\t@gsjj
   \def\@@{p}\ifx\@\@@\else
   \expandafter\gdef\csname#1cite\endcsname##1{\citejj{##1}}
   \expandafter\gdef\csname#1page\endcsname##1{?}
   \expandafter\gdef\csname#1tag\endcsname##1{\tagjj{##1}}\fi\fi\next}
\newif\ifshowstuffinmargin
\showstuffinmarginfalse
\def\jjtags{\showstuffinmargintrue
\ifx\all\relax\else\expandafter\t@gsjj\all,\fi}

\def\tagjj#1{\realstag{#1}\mginpar{\zeigen{#1}}}
\def\citejj#1{\zeigen{#1}\mginpar{\rechnen{#1}}}

\def\rechnen#1{\expandafter\ifx\csname stag#1\endcsname\relax ??\else
                           \csname stag#1\endcsname\fi}

\newdimen\theight

\def\marginfont{\sevenrm}

\def\trymarginbox#1{\setbox0=\hbox{\marginfont\hskip\marginshift #1}%
		\global\marginshift\wd0 
		\global\advance\marginshift\margindelta}

\def \mginpar#1{%
\ifvmode\setbox0\hbox to \hsize{\hfill\rlap{\marginfont\quad#1}}%
\ht0 0cm
\dp0 0cm
\box0\vskip-\baselineskip
\else 
             \vadjust{\trymarginbox{#1}%
		\ifdim\marginshift>\marginmax \global\marginshift\marginmin
			\trymarginbox{#1}%
                \fi
             \theight=\ht0
             \advance\theight by \dp0    \advance\theight by \lineskip
             \kern -\theight \vbox to \theight{\rightline{\rlap{\box0}}%
\vss}}\fi}


\def\t@gsoff#1,{\def\@{#1}\ifx\@\empty\let\next=\relax\else\let\next=\t@gsoff
   \def\@@{p}\ifx\@\@@\else
   \expandafter\gdef\csname#1cite\endcsname##1{\zeigen{##1}}
   \expandafter\gdef\csname#1page\endcsname##1{?}
   \expandafter\gdef\csname#1tag\endcsname##1{\zeigen{##1}}\fi\fi\next}
\def\verbatimtags{\showstuffinmarginfalse
\ifx\all\relax\else\expandafter\t@gsoff\all,\fi}
\def\zeigen#1{\hbox{$\langle$}#1\hbox{$\rangle$}}

\def\(#1){\edef\dot@g{\ifmmode\ifinner(\hbox{\noexpand\etag{#1}})
   \else\noexpand\eqno(\hbox{\noexpand\etag{#1}})\fi
   \else(\noexpand\ecite{#1})\fi}\dot@g}

\newif\ifbr@ck
\def\eat#1{}
\def\[#1]{\br@cktrue[\br@cket#1'X]}
\def\br@cket#1'#2X{\def\temp{#2}\ifx\temp\empty\let\next\eat
   \else\let\next\br@cket\fi
   \ifbr@ck\br@ckfalse\br@ck@t#1,X\else\br@cktrue#1\fi\next#2X}
\def\br@ck@t#1,#2X{\def\temp{#2}\ifx\temp\empty\let\neext\eat
   \else\let\neext\br@ck@t\def\temp{,}\fi
   \def\teemp{#1}\ifx\teemp\empty\else\rcite{#1}\fi\temp\neext#2X}
\def\resetbr@cket{\gdef\[##1]{[\rtag{##1}]}}
\def\references{\resetbr@cket\newsection References\par}

\newtoks\symb@ls\newtoks\s@mb@ls\newtoks\p@gelist\n@wcount\ftn@mber
    \ftn@mber=1\newif\ifftn@mbers\ftn@mbersfalse\newif\ifbyp@ge\byp@gefalse
\def\defm@rk{\ifftn@mbers\n@mberm@rk\else\symb@lm@rk\fi}
\def\n@mberm@rk{\xdef\m@rk{{\the\ftn@mber}}%
    \global\advance\ftn@mber by 1 }
\def\rot@te#1{\let\temp=#1\global#1=\expandafter\r@t@te\the\temp,X}
\def\r@t@te#1,#2X{{#2#1}\xdef\m@rk{{#1}}}
\def\b@@st#1{{$^{#1}$}}\def\str@p#1{#1}
\def\symb@lm@rk{\ifbyp@ge\rot@te\p@gelist\ifnum\expandafter\str@p\m@rk=1 
    \s@mb@ls=\symb@ls\fi\write\f@nsout{\number\count0}\fi \rot@te\s@mb@ls}
\def\byp@ge{\byp@getrue\n@wwrite\f@nsin\openin\f@nsin=\jobname.fns 
    \n@wcount\currentp@ge\currentp@ge=0\p@gelist={0}
    \re@dfns\closein\f@nsin\rot@te\p@gelist
    \n@wread\f@nsout\openout\f@nsout=\jobname.fns }
\def\m@kelist#1X#2{{#1,#2}}
\def\re@dfns{\ifeof\f@nsin\let\next=\relax\else\read\f@nsin to \f@nline
    \ifx\f@nline\v@idline\else\let\t@mplist=\p@gelist
    \ifnum\currentp@ge=\f@nline
    \global\p@gelist=\expandafter\m@kelist\the\t@mplistX0
    \else\currentp@ge=\f@nline
    \global\p@gelist=\expandafter\m@kelist\the\t@mplistX1\fi\fi
    \let\next=\re@dfns\fi\next}
\def\symbols#1{\symb@ls={#1}\s@mb@ls=\symb@ls} 
\def\bigsymbol{\textstyle}
\symbols{\bigsymbol\ast,\dagger,\ddagger,\sharp,\flat,\natural,\star}
\def\ftnumbers{\ftn@mberstrue} \def\ftsymbols{\ftn@mbersfalse}
\def\paginal{\byp@ge} \def\resetftnumbers{\ftn@mber=1}
\def\ftnote#1{\defm@rk\expandafter\expandafter\expandafter\footnote
    \expandafter\b@@st\m@rk{#1}}

\long\def\jump#1\endjump{}
\def\ssum{\mathop{\lower .1em\hbox{$\textstyle\Sigma$}}\nolimits}

\def\qed{\nobreak\kern 1em \vrule height .5em width .5em depth 0em}
\def\newneq{\hbox{\rlap{\hbox to 1\wd9{\hss$=$\hss}}\raise .1em 
   \hbox to 1\wd9{\hss$\scriptscriptstyle/$\hss}}}
\def\subsetne{\setbox9 = \hbox{$\subset$}\mathrel{\hbox{\rlap
   {\lower .4em \newneq}\raise .13em \hbox{$\subset$}}}}
\def\supsetne{\setbox9 = \hbox{$\subset$}\mathrel{\hbox{\rlap
   {\lower .4em \newneq}\raise .13em \hbox{$\supset$}}}}

\def\vbar{\mathchoice{\vrule height6.3ptdepth-.5ptwidth.8pt\kern-.8pt}
   {\vrule height6.3ptdepth-.5ptwidth.8pt\kern-.8pt}
   {\vrule height4.1ptdepth-.35ptwidth.6pt\kern-.6pt}
   {\vrule height3.1ptdepth-.25ptwidth.5pt\kern-.5pt}}
\def\f@dge{\mathchoice{}{}{\mkern.5mu}{\mkern.8mu}}
\def\b@c#1#2{{\rm \mkern#2mu\vbar\mkern-#2mu#1}}
\def\b@b#1{{\rm I\mkern-3.5mu #1}}
\def\b@a#1#2{{\rm #1\mkern-#2mu\f@dge #1}}
\def\bb#1{{\count4=`#1 \advance\count4by-64 \ifcase\count4\or\b@a A{11.5}\or
   \b@b B\or\b@c C{5}\or\b@b D\or\b@b E\or\b@b F \or\b@c G{5}\or\b@b H\or
   \b@b I\or\b@c J{3}\or\b@b K\or\b@b L \or\b@b M\or\b@b N\or\b@c O{5} \or
   \b@b P\or\b@c Q{5}\or\b@b R\or\b@a S{8}\or\b@a T{10.5}\or\b@c U{5}\or
   \b@a V{12}\or\b@a W{16.5}\or\b@a X{11}\or\b@a Y{11.7}\or\b@a Z{7.5}\fi}}

\catcode`\X=11 \catcode`\@=12


\expandafter\ifx\csname citeadd.tex\endcsname\relax
\expandafter\gdef\csname citeadd.tex\endcsname{}
\else \message{Hey!  Apparently you were trying to
\string\input{citeadd.tex} twice.   This does not make sense.} 
\errmessage{Please edit your file (probably \jobname.tex) and remove
any duplicate ``\string\input'' lines}\endinput\fi

\def\scitet{\sciteerror{ambiguous}}

\def\sciteerror#1#2{{\mathortextbf{\scite{#2}}}\complainaboutcitation{#1}{#2}}
\def\mathortextbf#1{\hbox{\bf #1}}
\def\complainaboutcitation#1#2{%
\vadjust{\line{\llap{---$\!\!>$ }\qquad scite$\{$#2$\}$ #1\hfil}}}

\sectno=-1   
\localtags
\NoBlackBoxes
\define\mr{\medskip\roster}
\define\sn{\smallskip\noindent}
\define\mn{\medskip\noindent}
\define\bn{\bigskip\noindent}
\define\ub{\underbar}
\define\wilog{\text{without loss of generality}}
\define\ermn{\endroster\medskip\noindent}
\define\dbca{\dsize\bigcap}
\define\dbcu{\dsize\bigcup}
\define \nl{\newline}
\magnification=\magstep 1
\documentstyle{amsppt}

{    
\catcode`@11

\ifx\alicetwothousandloaded@\relax
  \endinput\else\global\let\alicetwothousandloaded@\relax\fi

\gdef\subjclass{\let\savedef@\subjclass
 \def\subjclass##1\endsubjclass{\let\subjclass\savedef@
   \toks@{\def\usualspace{{\rm\enspace}}\eightpoint}%
   \toks@@{##1\unskip.}%
   \edef\thesubjclass@{\the\toks@
     \frills@{{\noexpand\rm2000 {\noexpand\it Mathematics Subject
       Classification}.\noexpand\enspace}}%
     \the\toks@@}}%
  \nofrillscheck\subjclass}
} 

\pageheight{8.5truein}
\topmatter
\title{Superatomic Boolean Algebras: maximal rigidity} \endtitle
\rightheadtext{Moving atoms}
\author {Saharon Shelah \thanks {\null\newline I would like to thank 
Alice Leonhardt for the beautiful typing. \null\newline
 Latest Revision - 00/July/11 \null\newline
 sh704 paper} \endthanks} \endauthor 
\affil{Institute of Mathematics\\
 The Hebrew University\\
 Jerusalem, Israel
 \medskip
 Rutgers University\\
 Mathematics Department\\
 New Brunswick, NJ  USA} \endaffil

\keywords  Set Theory, Boolean Algebras, superatomic, rigid; pcf, MAD \endkeywords

\abstract   We prove that for any superatomic Boolean Algebra of cardinality
$> \beth_\omega$ there is an automorphism moving uncountably many atoms.
Similarly for larger cardinals any of those results are essentially best
possible  \endabstract
\endtopmatter
\document  

\expandafter\ifx\csname alice2jlem.tex\endcsname\relax
  \expandafter\xdef\csname alice2jlem.tex\endcsname{\the\catcode`@}
\else \message{Hey!  Apparently you were trying to
\string\input{alice2jlem.tex}  twice.   This does not make sense.}
\errmessage{Please edit your file (probably \jobname.tex) and remove
any duplicate ``\string\input'' lines}\endinput\fi

\expandafter\ifx\csname bib4plain.tex\endcsname\relax
  \expandafter\gdef\csname bib4plain.tex\endcsname{}
\else \message{Hey!  Apparently you were trying to \string\input
  bib4plain.tex twice.   This does not make sense.}
\errmessage{Please edit your file (probably \jobname.tex) and remove
any duplicate ``\string\input'' lines}\endinput\fi

\def\renewcommand{\newcommand}	       
\edef\cite{\the\catcode`@}%
\catcode`@ = 11
\let\@oldatcatcode = \cite
\chardef\@letter = 11
\chardef\@other = 12
%
%
%
%
\def\@innerdef#1#2{\edef#1{\expandafter\noexpand\csname #2\endcsname}}%
%
%
\@innerdef\@innernewcount{newcount}%
\@innerdef\@innernewdimen{newdimen}%
\@innerdef\@innernewif{newif}%
\@innerdef\@innernewwrite{newwrite}%
%
%
%
\def\@gobble#1{}%
%
%
%
\ifx\inputlineno\@undefined
   \let\@linenumber = \empty 
\else
   \def\@linenumber{\the\inputlineno:\space}%
\fi
%
%
%
\def\@futurenonspacelet#1{\def\cs{#1}%
   \afterassignment\@stepone\let\@nexttoken=
}%
\begingroup 
\def\\{\global\let\@stoken= }%
\\ 
\endgroup
\def\@stepone{\expandafter\futurelet\cs\@steptwo}%
\def\@steptwo{\expandafter\ifx\cs\@stoken\let\@@next=\@stepthree
   \else\let\@@next=\@nexttoken\fi \@@next}%
\def\@stepthree{\afterassignment\@stepone\let\@@next= }%
%
%
%
\def\@getoptionalarg#1{%
   \let\@optionaltemp = #1%
   \let\@optionalnext = \relax
   \@futurenonspacelet\@optionalnext\@bracketcheck
}%
%
%
\def\@bracketcheck{%
   \ifx [\@optionalnext
      \expandafter\@@getoptionalarg
   \else
      \let\@optionalarg = \empty
      \expandafter\@optionaltemp
   \fi
}%
\def\@@getoptionalarg[#1]{%
   \def\@optionalarg{#1}%
   \@optionaltemp
}%
%
%
%
\def\@nnil{\@nil}%
\def\@fornoop#1\@@#2#3{}%
\def\@for#1:=#2\do#3{%
   \edef\@fortmp{#2}%
   \ifx\@fortmp\empty \else
      \expandafter\@forloop#2,\@nil,\@nil\@@#1{#3}%
   \fi
}%
\def\@forloop#1,#2,#3\@@#4#5{\def#4{#1}\ifx #4\@nnil \else
       #5\def#4{#2}\ifx #4\@nnil \else#5\@iforloop #3\@@#4{#5}\fi\fi
}%
\def\@iforloop#1,#2\@@#3#4{\def#3{#1}\ifx #3\@nnil
       \let\@nextwhile=\@fornoop \else
      #4\relax\let\@nextwhile=\@iforloop\fi\@nextwhile#2\@@#3{#4}%
}%
%
%
%
\@innernewif\if@fileexists
\def\@testfileexistence{\@getoptionalarg\@finishtestfileexistence}%
\def\@finishtestfileexistence#1{%
   \begingroup
      \def\extension{#1}%
      \immediate\openin0 =
         \ifx\@optionalarg\empty\jobname\else\@optionalarg\fi
         \ifx\extension\empty \else .#1\fi
         \space
      \ifeof 0
         \global\@fileexistsfalse
      \else
         \global\@fileexiststrue
      \fi
      \immediate\closein0
   \endgroup
}%
%
%
%
%
\def\bibliographystyle#1{%
   \@readauxfile
   \@writeaux{\string\bibstyle{#1}}%
}%
\let\bibstyle = \@gobble
%
%
\let\bblfilebasename = \jobname
\def\bibliography#1{%
   \@readauxfile
   \@writeaux{\string\bibdata{#1}}%
   \@testfileexistence[\bblfilebasename]{bbl}%
   \if@fileexists
      \nobreak
      \@readbblfile
   \fi
}%
\let\bibdata = \@gobble
%
%
\def\nocite#1{%
   \@readauxfile
   \@writeaux{\string\citation{#1}}%
}%
\@innernewif\if@notfirstcitation
%
%
\def\cite{\@getoptionalarg\@cite}%
%
%
\def\@cite#1{%
   \let\@citenotetext = \@optionalarg
   \printcitestart
   \nocite{#1}%
   \@notfirstcitationfalse
   \@for \@citation :=#1\do
   {%
      \expandafter\@onecitation\@citation\@@
   }%
   \ifx\empty\@citenotetext\else
      \printcitenote{\@citenotetext}%
   \fi
   \printcitefinish
}%
\def\@onecitation#1\@@{%
   \if@notfirstcitation
      \printbetweencitations
   \fi
   \expandafter \ifx \csname\@citelabel{#1}\endcsname \relax
      \if@citewarning
         \message{\@linenumber Undefined citation `#1'.}%
      \fi
      \expandafter\gdef\csname\@citelabel{#1}\endcsname{%
\strut
\vadjust{\vskip-\dp\strutbox
\vbox to 0pt{\vss\parindent0cm \leftskip=\hsize 
\advance\leftskip3mm
\advance\hsize 4cm\strut\openup-4pt 
\rightskip 0cm plus 1cm minus 0.5cm ?  #1 ?\strut}}
         {\tt
            \escapechar = -1
            \nobreak\hskip0pt
            \expandafter\string\csname#1\endcsname
            \nobreak\hskip0pt
         }%
      }%
   \fi
   \csname\@citelabel{#1}\endcsname
   \@notfirstcitationtrue
}%
%
%
\def\@citelabel#1{b@#1}%
%
%
\def\@citedef#1#2{\expandafter\gdef\csname\@citelabel{#1}\endcsname{#2}}%
%
%
%
\def\@readbblfile{%
   \ifx\@itemnum\@undefined
      \@innernewcount\@itemnum
   \fi
   \begingroup
      \def\begin##1##2{%
         \setbox0 = \hbox{\biblabelcontents{##2}}%
         \biblabelwidth = \wd0
      }%
      \def\end##1{}
      %
      %
      \@itemnum = 0
      \def\bibitem{\@getoptionalarg\@bibitem}%
      \def\@bibitem{%
         \ifx\@optionalarg\empty
            \expandafter\@numberedbibitem
         \else
            \expandafter\@alphabibitem
         \fi
      }%
      \def\@alphabibitem##1{%
         \expandafter \xdef\csname\@citelabel{##1}\endcsname {\@optionalarg}%
         \ifx\biblabelprecontents\@undefined
            \let\biblabelprecontents = \relax
         \fi
         \ifx\biblabelpostcontents\@undefined
            \let\biblabelpostcontents = \hss
         \fi
         \@finishbibitem{##1}%
      }%
      \def\@numberedbibitem##1{%
         \advance\@itemnum by 1
         \expandafter \xdef\csname\@citelabel{##1}\endcsname{\number\@itemnum}%
         \ifx\biblabelprecontents\@undefined
            \let\biblabelprecontents = \hss
         \fi
         \ifx\biblabelpostcontents\@undefined
            \let\biblabelpostcontents = \relax
         \fi
         \@finishbibitem{##1}%
      }%
      \def\@finishbibitem##1{%
         \biblabelprint{\csname\@citelabel{##1}\endcsname}%
         \@writeaux{\string\@citedef{##1}{\csname\@citelabel{##1}\endcsname}}%
         \ignorespaces
      }%
      %
      %
      \let\em = \bblem
      \let\newblock = \bblnewblock
      \let\sc = \bblsc
      \frenchspacing
      \clubpenalty = 4000 \widowpenalty = 4000
      \tolerance = 10000 \hfuzz = .5pt
      \everypar = {\hangindent = \biblabelwidth
                      \advance\hangindent by \biblabelextraspace}%
      \bblrm
      \parskip = 1.5ex plus .5ex minus .5ex
      \biblabelextraspace = .5em
      \bblhook
      \input \bblfilebasename.bbl
   \endgroup
}%
%
%
\@innernewdimen\biblabelwidth
\@innernewdimen\biblabelextraspace
%
%
%
\def\biblabelprint#1{%
   \noindent
   \hbox to \biblabelwidth{%
      \biblabelprecontents
      \biblabelcontents{#1}%
      \biblabelpostcontents
   }%
   \kern\biblabelextraspace
}%
%
%
%
\def\biblabelcontents#1{{\bblrm [#1]}}%
%
%
\def\bblrm{\rm}%
%
%
\def\bblem{\it}%
%
%
\def\bblsc{\ifx\@scfont\@undefined
              \font\@scfont = cmcsc10
           \fi
           \@scfont
}%
%
%
\def\bblnewblock{\hskip .11em plus .33em minus .07em }%
%
%
\let\bblhook = \empty
%
%
%
\def\printcitestart{[}
\def\printcitefinish{]}
\def\printbetweencitations{, }
\def\printcitenote#1{, #1}
%
%
%
\let\citation = \@gobble
%
%
%
\@innernewcount\@numparams
%
%
\def\newcommand#1{%
   \def\@commandname{#1}%
   \@getoptionalarg\@continuenewcommand
}%
%
%
\def\@continuenewcommand{%
   \@numparams = \ifx\@optionalarg\empty 0\else\@optionalarg \fi \relax
   \@newcommand
}%
%
%
\def\@newcommand#1{%
   \def\@startdef{\expandafter\edef\@commandname}%
   \ifnum\@numparams=0
      \let\@paramdef = \empty
   \else
      \ifnum\@numparams>9
         \errmessage{\the\@numparams\space is too many parameters}%
      \else
         \ifnum\@numparams<0
            \errmessage{\the\@numparams\space is too few parameters}%
         \else
            \edef\@paramdef{%
               \ifcase\@numparams
                  \empty  No arguments.
               \or ####1%
               \or ####1####2%
               \or ####1####2####3%
               \or ####1####2####3####4%
               \or ####1####2####3####4####5%
               \or ####1####2####3####4####5####6%
               \or ####1####2####3####4####5####6####7%
               \or ####1####2####3####4####5####6####7####8%
               \or ####1####2####3####4####5####6####7####8####9%
               \fi
            }%
         \fi
      \fi
   \fi
   \expandafter\@startdef\@paramdef{#1}%
}%
%
%
%
%
\def\@readauxfile{%
   \if@auxfiledone \else 
      \global\@auxfiledonetrue
      \@testfileexistence{aux}%
      \if@fileexists
         \begingroup
            \endlinechar = -1
            \catcode`@ = 11
            \input \jobname.aux
         \endgroup
      \else
         \message{\@undefinedmessage}%
         \global\@citewarningfalse
      \fi
      \immediate\openout\@auxfile = \jobname.aux
   \fi
}%
%
%
\newif\if@auxfiledone
\ifx\noauxfile\@undefined \else \@auxfiledonetrue\fi
%
%
%
%
\@innernewwrite\@auxfile
\def\@writeaux#1{\ifx\noauxfile\@undefined \write\@auxfile{#1}\fi}%
%
%
%
\ifx\@undefinedmessage\@undefined
   \def\@undefinedmessage{No .aux file; I won't give you warnings about
                          undefined citations.}%
\fi
%
%
\@innernewif\if@citewarning
\ifx\noauxfile\@undefined \@citewarningtrue\fi
%
%
%
\catcode`@ = \@oldatcatcode


\def\widestnumber#1#2{}

\def\rm{\fam0 \tenrm}

\def\fakesubhead#1\endsubhead{\bigskip\noindent{\bf#1}\par}



%
%
%

%

\font\textrsfs=rsfs10
\font\scriptrsfs=rsfs7
\font\scriptscriptrsfs=rsfs5

\newfam\rsfsfam
\textfont\rsfsfam=\textrsfs
\scriptfont\rsfsfam=\scriptrsfs
\scriptscriptfont\rsfsfam=\scriptscriptrsfs

\edef\oldcatcodeofat{\the\catcode`\@}
\catcode`\@11

\def\Cal@@#1{\noaccents@ \fam \rsfsfam #1}

\catcode`\@\oldcatcodeofat


\expandafter\ifx \csname margininit\endcsname \relax\else\margininit\fi

\mn

\head {Anotated Content} \endhead  \resetall 
\bn
\S1  $\quad$ Superatomic Boolean Algebra have nontrivial automorphism
\mr
\item "{{}}"  [We prove that if $B$ is a superatomic Boolean Algebra, then
it has quite a nontrivial automorphism; specifically if $B$ is of cardinality
$> \beth_4(\sigma)$ then $B$ has an automorphism moving $> \sigma$
atoms.  We then discuss how much we can weaken the superatomicity assumptions.]
\ermn
\S2 $\quad$  Constructing counterexamples
\mr
\item "{{}}" [Under some assumptions we construct examples of superatomic
Boolean Algebras for which every automorphism move few atoms.]
\ermn
\S3 $\quad$  Sufficient conditions for the existence of $\langle B_i:i <
\mu \rangle$
\mr
\item "{{}}"  [We deal with the assumptions of the construction in \S2.]
\endroster
\newpage

\head {\S0 Introduction} \endhead  \resetall \sectno=0
\bigskip

We show that for a superatomic Boolean Algebra has an automorphism
moving uncountably many atoms if it is large enough, really $>
\beth_u$; similarly replacing $\aleph_0$ by $\theta$; (an automorphism
move an atom if its image is not itself).  We then show that those
results are essentially best possible.  Of course, we can express
those results in topological terms.  See \cite{M} and his cite on
background and history, in particular work of Rubin and work of the
author.
\bigskip

\noindent {\bf Notation} 
\bigskip

\definition{\stag{0.1} Definition}  1) For a Boolean Algebra $B$ let us
define the ideal id$_\alpha(B)$ by induction:
\mr
\item "{${{}}$}"  id$_0(B) = \{0\}$
\sn
\item "{${{}}$}"  id$_\beta(B) = \{x_1 \cup \ldots \cup x_n:\text{ for
some } \alpha < \beta \text{ and } n < \omega \text{ for each } \ell \in
\{1,\dotsc,n\}$ the element $x_\ell/\text{id}_\alpha(B)$ is an atom of 
$B/\text{id}_\alpha(B)$ or $x_\ell \in \text{ id}_\alpha(B)\}$.
\ermn
Hence for limit $\delta$ we have
\mr
\item "{${{}}$}"  id$_\delta(B) = \dbcu_{\beta < \delta} \text{ id}_\beta(B)$
for limit $\delta$. \nl
Let id$_\infty(B) = \dbcu_\alpha \text{ id}_\alpha(B)$.
\ermn
2) For $x \in \text{ id}_\infty(B)$ let rk$(x,B) = \text{ Min}\{\alpha:
x \in \text{ id}_{\alpha +1}(B)\}$. \nl
3) $B$ is superatomic if $B = \text{ id}_\infty(B)$ and dp$(B)$ be the
ordinal $\alpha$ such that $B/\text{id}_\alpha(B)$ is a finite Boolean Algebra.
\enddefinition
\newpage

\head {\S1 Superatomic Boolean Algebra's have nontrivial automorphisms}
\endhead  \resetall 
\bigskip

\proclaim{\stag{gd.1} Theorem}  Assume
\mr
\item "{$(a)$}"  $B$ is a superatomic Boolean Algebra with no automorphism
moving $\ge \theta$ atoms; that is if $\pi$ ia an automorphism of $B$ then
\nl
$\{x:x \in \text{ atom}(B) \text{ and } \pi(x) \ne x\}$ is a set of
cardinality $< \theta$
\sn
\item "{$(b)$}"  $\theta$ regular uncountable.
\ermn
\ub{Then}  $|B| \le \beth_3(2^{< \theta})$, so if $\theta = \sigma^+$ then
$|B| \le \beth_4(\sigma)$.
\endproclaim
\bigskip

\demo{Proof}  Toward contradiction let $B$ be a counterexample and let
$\mu$ be the number of atoms of $B$.
Without loss of generality
\mr
\item "{$\boxtimes_1$}"  $B$ is a Boolean Algebra of subsets of $\mu$ and
its atoms are the singletons $\{\alpha\},\alpha < \mu$.
\ermn
Let $I =: [\mu]^{< \theta} \cap B = \{x \in B:|x| < \theta\}$, clearly $I$ is
an ideal of $B$ and let \nl
$Y =: \{x:x \in B \text{ and } x/I \text{ is an
atom of } B/I\}$. \nl
We shall prove (after some preliminary things) that:
\mr
\item "{$\boxtimes_2$}"   if $x \in Y$ then $|x| \le 2^{(2^{<\theta})}$.
\ermn
We shall say that a set $a \subseteq \mu$ is $B$-\ub{autonomous} if 
$(\forall y \in I)(y \cap a \in B)$; in this case we let 
$B \restriction a = B \cap {\Cal P}(a)$.
\nl
Clearly
\mr
\item "{$\otimes_1$}"  the family of $B$-autonomous subsets of $\mu$ is a
Boolean ring of subsets of $\mu$ (i.e. closed under $a \cap b,a \cup b,
a \backslash b$) and include $I$
\sn
\item "{$\otimes_2$}"  for a $B$-autonomous $a,B \restriction a$ is a 
Boolean algebra of subsets of $a$ which include $\{\{\alpha\}:\alpha \in a\}$.
\ermn
Also
\mr
\item "{$\otimes_3$}"  if $a_0,a_1$ are $B$-autonomous subsets of $\mu,
x \in Y,a_0 \subseteq x,a_1 \subseteq x$ and $B \restriction a_0 \cong B
\restriction a_1$ over $B \restriction (a_1 \cap a_2) = B \cap 
{\Cal P}(a_1 \cap a_2)$, \ub{then} there is 
an automorphism $h$ of $B$ such that
$h$ maps $a_0$ to $a_1,a_1$ to $a_0$ and $\alpha \in \mu \backslash a_0
\backslash a_1 \Rightarrow h(\{\alpha\}) = \{\alpha\}$. \nl
\sn
[Why?  Let $g$ be an isomorphism from $B \restriction a_0$ 
onto $B \restriction a_1$ over $B \restriction (a_0 \cap a_1)$; now 
we define a permutation $h$ of
atom$(B) = \{\{\alpha\}:\alpha < \mu\}$; let $\alpha \in a_0
\Rightarrow h(\{\alpha\}) = g(\{\alpha\}),h(g(\{\alpha\})) = \{\alpha\}$ and
$\alpha \in \mu \backslash a_0 \backslash a_1 \Rightarrow h(\{\alpha\}) =
\{\alpha\}$, by the demands on $g$ clearly $h$ is a well defined permutation
of atom$(B)$.  Now $h$ can be naturally extended to an automorphism $\hat h$ 
of ${\Cal P}(\mu)$ of order two, we have to check that $\hat h$ maps $B$ onto
itself; even into itself suffice (because of ``order two").  Clearly $\hat h(x) = x$ and
$\hat h \restriction (B \restriction (\mu \backslash x))$ is the identity.  
So it is enough to check : that 
$\hat h \restriction (B \restriction x)$ is an automorphism of $B \restriction
x$.  But $I \cap (B \restriction x)$ is a maximal ideal of $B \restriction x$ 
(as $x \in Y$) hence it is enough to check that $\hat h$ maps 
$I \cap (B \restriction x)$
into itself.  As $b \in I \cap (B \restriction x) \Rightarrow
b = (b \backslash a_0 \backslash a_1) \cup (b \cap a_0 \cap a_1) 
\cup (b \cap a_0 \backslash a_1) \cup (b_1 \cap a_1 \backslash a_0)$, and all
four are in $I$; \ub{clearly} it is enough to check the following
statements: $b \in I \and b 
\subseteq x \backslash a_0 \backslash a_1 \Rightarrow h(b) \in I$, and 
$\ell < 2 \and b \in I \and
b \subseteq x \cap a_\ell \backslash a_{1 - \ell} \Rightarrow \hat h(b) \in
I$ and lastly $b \in I \and b \subseteq a_0 \cap a_1 \Rightarrow \hat h(b) 
\in I$.  The second implication holds by the choice of $g$, the 
first as $\hat h(b) = b$ in this case and the last one as $h \restriction
(a_0 \cap a_1)$ is the identity so again $\hat h(b) = b$.]
\sn
\item "{$\otimes_4$}"  if $b \subseteq \mu,|b| \le 2^{< \theta}$ then for
some $B$-autonomous set $c$ we have $a \subseteq c \subseteq \mu,|c| \le
2^{< \theta}$. \nl
[Why?  Find $c$ satisfying $b \subseteq c \subseteq \mu,
|c| \le 2^{< \theta}$ such that
$(\forall y \in [c]^{< \theta})[(\exists z)(y \subseteq z \in I) \rightarrow
(\exists z \subseteq c)(y \subseteq z \in I)]$, just close $\theta$ times
recalling $\theta$ is regular.  Now if $g \in I$ then $|g| < \theta$
hence $g \cap c \in [c]^{< \theta}$ so there is $z$ such that $g \cap
c \in z \in I \and z \subseteq c$; hence $y \cap c = y \cap z \in I$.
This proves that $c$ is autonomous.  Now check.]
\endroster
\enddemo
\bn
Now we return to the promised $\boxtimes_2$.
\demo{Proof of $\boxtimes_2$}  if 
$|x| > 2^{2^{< \theta}}$ let $\alpha_i \in x$
for $i < (2^{(2^{< \theta})})^+$ be pairwise distinct, let $a_i$ be
$B$-autonomous set of cardinality $\le 2^{< \theta}$ such that
$\{\alpha_{i + \varepsilon}:\varepsilon < 2^{< \theta}\} \subseteq a_i$
(exists by $\otimes_4$), and \wilog \, $a_i \subseteq x$ (just use $a_i \cap
x$).  For some club $C$ of $(2^{2^{< \theta}})^+$, we have $i < j \in C
\Rightarrow a_i \cap \{\alpha_{j + \varepsilon}:\varepsilon < 2^{< \theta}\}
= \emptyset$ hence $i < j \in C \Rightarrow |a_j \backslash a_i| \ge
2^{< \theta}$.
Now $I \cap {\Cal P}(a_i)$ has cardinality $\le |a_i|^{< \theta} \le
2^{< \theta}$ (as $\theta$ is regular) but $x \in Y$ hence 
$B \restriction a_i$ has
cardinality $\le 2^{< \theta}$,  hence the number of isomorphism types of
$B \restriction a_i$ is $\le 2^{(2^{< \theta})}$.  Hence there is a
stationary $S \subseteq \{\delta < (2^{(2^{<
\theta})})^+:\text{cf}(\delta) = (2^{< \theta})^+\}$ and $a^*$ such
that $i \in S \and j \in S \and i \ne j \Rightarrow a_i \cap a_j =
a^*$ (the $\Delta$-system lemma).  Also the number of isomorphic types
of $(B \restriction a_i,\alpha)_{\alpha \in a^*}$ is as most $\le
2^{(2^{< \theta})}$ hence for some $i < j$ from $C \cap S$ 
we have $B \restriction a_i \cong B \restriction a_j$, but $|a_j
\backslash a_i| \ge 2^{< \theta} \ge \theta$ hence by $\otimes_3$ there is an
automorphism $h$ of $B$ which moves $\ge 2^{< \theta}$ atoms, contradiction. 
\nl
Next
\mr
\item "{$\boxtimes_3$}"  $|Y/I| \le \beth_2(2^{< \theta})$. \nl
[Why?  If not, we can find $x_i \in Y$ for $i < (\beth_2(2^{< \theta}))^+$ 
such that $i \ne j \Rightarrow x_i/I \ne x_j/I$.
As $|x_i| \le \beth_1(2^{< \theta})$ by $\boxtimes_2$, by the $\Delta$-system
lemma for some unbounded $A \subseteq (\beth_2(2^{< \theta}))^+$ the set
$\langle x_i:i \in A \rangle$ is a $\Delta$-system hence \wilog \, $\langle x_i:i \in A
\rangle$ are pairwise disjoint (not really needed just clearer).  As $B
\restriction x_i$ is a Boolean Algebra of cardinality $\le \beth_1
(2^{< \theta})$ (as $I \cap {\Cal P}(x_i)$ is a maximal ideal of $B
\restriction x_i$ and $I \cap {\Cal P}(x_i) \subseteq [x_i]^{< \theta}$ and
$|x_i| \le \beth_1(2^{< \theta})$ by $\boxtimes_2$)
there are at most $\beth_2(2^{< \theta})$ isomorphism types of $B
\restriction x_i$.  So for some $i \ne j$ in $A$ we have $B \restriction x_i
\cong B \restriction x_j$, so as in the proof of $\otimes_3$ there is an
automorphism $h$ of $B$ mapping $x_i$ to $x_j$ hence moving $\ge |x_i
\backslash x_j| \ge \theta$ atoms because $x_i \ne x_j$ mod $I$.]
\ermn
Choose a set $\{x_\alpha:\alpha < \alpha^* \le \beth_2(2^{< \theta})\}$ of
representatives of $Y/I$ and let $x^* = \dbcu_{\alpha < \alpha^*} x_\alpha$,
so $x^* \subseteq \mu,|x^*| \le \beth_2(2^{< \theta})$. \nl
Define $J = \{a \in B:a \cap x^* = \emptyset\}$.
\mr
\item "{$\boxtimes_4$}"  $J \subseteq I$. \nl
[Why?  If not, there is $x \in J \backslash I$ such that $x/I$ is an atom
of $B/I$ so $x/I \in \{x_\alpha/I:\alpha < \alpha^*\}$, so for some $\alpha,
x/I = x_\alpha/I$ hence $|x \backslash x_\alpha| < \theta$ hence $|x \cap
x_\alpha| = \theta$ hence $x \cap x^* \ne \emptyset$ hence $x \notin J$,
a contradiction.]
\ermn
Define an equivalent relation ${\Cal E}$ on $B:y_1 {\Cal E} y_2$ iff
$y_1 \cap x^* = y_2 \cap x^*$.  Clearly ${\Cal E}$ has $\le 2^{|x^*|}$
equivalence classes and $2^{|x^*|} \le \beth_3(2^{< \theta})$; also $y_1
{\Cal E} y_2 \rightarrow y_1 \backslash y_2 \in J$ (see its definition).  
Choose a set of representatives $\{y_\gamma:\gamma < \gamma^*\}$ for 
${\Cal E}$ so $\gamma^* \le \beth_3(2^{< \theta})$ and let $B^*$ be the
subalgebra of $B$ which $\{y_\gamma:\gamma < \gamma^*\}$ generates.  So
$|B^*| \le \beth_3(2^{< \theta})$ and, being superatomic, the number of 
ultrafilters of $B^*$ is also $\le \beth_3(2^{< \theta})$.  Next 
$B$ is generated by $J \cup B^*$ as for $y \in B$ there is $\gamma$ such that
$y {\Cal E} y_\gamma$ and $y_\gamma \in B^*,y - y_\gamma \in J,y_\gamma -
y \in J$ hence $y \in \langle J \cup B^* \rangle$.  For
$D$ an ultrafilter of $B^*$ let $Z_D = \{\alpha \in \mu:(\forall y \in B^*)
(\alpha \in y \leftrightarrow y \in D)\}$. \nl
Clearly
\mr
\item "{$\boxtimes_5$}"  for every $\alpha \in \mu \backslash x^*$ there is
a unique ultrafilter $D = D[\alpha]$ on $B^*$ such that $\alpha \in
Z_D$ (and the number of such ultrafilters is $\le \beth_3 (2^{<
\theta}))$.
\ermn
Now
\mr
\item "{$\boxtimes_6$}"  $\mu \le \beth_3(2^{< \theta})$. \nl
[Why?  Assume that not.  By $\otimes_4$ 
for each $i < \mu$ we can find a $B$-autonomous $a_i$ such
that $|a_i| \le 2^{< \theta}$ and $[i,i+2^{< \theta}) \subseteq a_i$; let
$a_i = \{\beta_{i,\varepsilon}:\varepsilon < \varepsilon_i\}$ with
$\beta_{i,\varepsilon}$ increasing with
$\varepsilon$.  Without loss of generality for some unbounded $A \subseteq
(\beth_3(2^{< \theta}))^+$ for all $i \in A$ the following does not depend
on $i:\varepsilon_i$ and $D[\beta_{i,\varepsilon}]$ for $\varepsilon <
\varepsilon_i$ (use $\boxtimes_5$), 
and $\{u \in [\varepsilon_i]^{< \theta}:\{\beta_{\varepsilon,
i}:i \in u\} \in I\}$, the $\varepsilon$ such that $\beta_{i,\varepsilon}=i$
and \wilog \, for $j < i$ in $A,a_j \cap [i,i + 2^{< \theta}) = \emptyset$.
By the $\Delta$-system lemma \wilog \, for some $a^*$ we have:
for $i < j$ in $A,a_i \cap a_j = a^*$.  So by $\otimes_1$ the set $a^*$ is
$B$-autonomous and also $a_i \backslash a^*$ is 
so we can use $a_i \backslash a^*$, so \wilog \, for $i \ne j$
in $A,a_i \cap a_j = \emptyset$ and as $|x^*| \le \beth_3 (2^{<
\theta})$ clearly \wilog \, $i \in A \Rightarrow a_i \cap
x^* = \emptyset$.
So for $i \ne j$ in $A$ there is an automorphism of $B$ interchanging $a_i,
a_j$ (the proof is like that proof of $\otimes_3$ using ``$B$ is
generated by $J \cup B^*$").  So we get a contradiction.]
\ermn
So $|J| \le |[\mu]^{< \theta}| = \mu^{< \theta} \le (\beth_3(2^{< \theta}))
^{< \theta} = \beth_3(2^{< \theta})$ so as $B$ is generated by $J \cup B^*$
together we get the desired conclusion. \hfill$\square_{\scite{gd.1}}$
\enddemo
\bigskip

\demo{\stag{gd.1a} Discussion}  1) We can weaken the assumption ``$B$
is superatomic by $B/I_{< \theta}[B]$ is superatomic", where: for a
Boolean Algebra $B$ and infinite cardinal $\theta$ we define $I_{<
\theta}(B) = \{x \in B:B \restriction x$ has density $< \theta\}$ (see
a little in \cite[\S1]{Sh:397}).  For $B$ superatomic this is the $I$
in the proof of \scite{gd.1} and the proof is the same.  What if we
just assume ``$B/I_{< \theta}[B]$ is atomic"?  One point in the proof
may fail: the number of ultrafilters of $B^*$ is not $\le
|\text{atom}(B^*)| \le \beth_3(2^{< \theta})$ but is $\le 2^{|B^*|}
\le 2^{2^{|Y|}} \le \beth_4(2^{< \theta})$, so we should replace
$\beth_3(2^{< \theta})$ by $\beth_4(2^{< \theta})$ in the conclusion.
We can adapt \scite{gd.2} to this case: e.g. let $\langle
d_\zeta:\zeta < \lambda = 2^\mu \rangle$ be a family of subsets of
$\mu$ such that any finite Boolean combination of them is infinite and
let $B^*$ be the Boolean subalgebra of ${\Cal P}[\mu]$ generated by
$\{c_\alpha:\alpha < \lambda = 2^\mu\} \cup \{\{i\}:i < \mu\}$. We let
$\lambda' = 2^\lambda$ let $\{c_\gamma:\gamma < \lambda'\}$ be an
independent family of subsets of $\lambda$ and we let $X^* =
\dbcu_{\alpha < \mu} X_\alpha \cup \{X^*_\gamma:\gamma < \lambda\}$?
We ignore $A'$ (and omit clause (k) of the assumption) and among the
generators of $\bold B$, clause (i), (ii) remains and
\mr
\item "{$(iii)'$}"  $c_\zeta = \{x \in X$: for some $\alpha \in
d_\zeta$ we have $x \in X_\alpha\} \cup \{X^*_\gamma:\zeta \in
C_\gamma,\gamma \in [\mu,\lambda')$.
\endroster
\enddemo
\newpage

\head {\S2 Constructing counterexamples} \endhead  \resetall \sectno=2
\bigskip

We would like to show that the bound from \S1 is essentially best possible.
The construction (in \scite{gd.2}) is closely related to the proof in
\S1, but we need various assumptions.  We shall deal with them later.
\bigskip

\proclaim{\stag{gd.2} Lemma}  Assume
\mr
\item "{$(a)$}"  $\theta \le \kappa \le \mu \le \lambda' \le \lambda,
\aleph_0 < \theta,\theta = {\text{\rm cf\/}}(\theta) = \sigma \ge \aleph_1$
\sn
\item "{$(b)$}"  there is an ${\Cal A} \subseteq [\mu]^{\aleph_0}$ almost
disjoint  (i.e. $A \ne B \in {\Cal A} \Rightarrow |A \cap B| < \aleph_0$)
such that $(\forall A \in [\mu]^\sigma)(\exists B \in {\Cal A})
(B \subseteq^* A)$ and $|{\Cal A}| = \mu$
\sn
\item "{$(c)$}"  $\bar B = \langle B_\alpha:\alpha < \mu \rangle$
\sn
\item "{$(d)$}"  $B_\alpha$ is a superatomic Boolean Algebra with $\le
\kappa$ atoms such that any automorphism of $B_\alpha$ moves $< \theta$ atoms
and $|B_\alpha| \le \lambda$; moreover if $c_1,c_2 \in I_\alpha$ (see
below) and $f$ is an isomorphism from $B_\alpha \restriction (1 -
c_1)$ onto $B_\alpha \restriction (1 - c_2)$ then $\theta > |\{x \in
\text{ atom}(B_\alpha):x \le_{B_\alpha} c_1$ or $f(x) \ne x\}|$
\sn
\item "{$(e)$}"   $I_\alpha = \{b \in B_\alpha:|\{x \in$ 
{\rm atom}$(B_\alpha):x \le b\}| < \theta\}$ is a maximal ideal of $B_\alpha$
\sn
\item "{$(f)$}"  there is an infinite set $X$ of atoms of $B_\alpha$ such
that for every $a \in B_\alpha,\{x \in X:x \le a\}$ is a finite or co-finite
subset of $X$
\sn
\item "{$(g)$}"  if $\alpha \ne \beta$ then for no $a_\alpha \in I_\alpha,
a_\beta \in I_\beta$ do we have \nl
$B_\alpha \restriction (1_{B_\alpha} - a_\alpha) \cong B_\beta \restriction
(1_{B_\beta} - a_\beta)$
\sn
\item "{$(h)$}"  $B^*$ is a superatomic Boolean Algebra
\sn
\item "{$(i)$}"  $B^*$ has $\mu$ atoms and $\lambda$ elements 
\footnote{if there is a tree ${\Cal T}$ with $\le \mu$
nodes and $\ge \lambda$ branches (= maximal linearly ordered subsets) then
such $B^*$ exists}
\sn
\item "{$(j)$}"  if ${\Cal U}$ is an infinite set of atoms of $B^*$ \ub{then}
for some $b \in B^*$ we have: [used?]
{\roster
\itemitem{ $(i)$ }  $\{x \in {\Cal U}:x \le b\}$ is infinite
\sn
\itemitem{ $(ii)$ }  $b/\text{id}_{\text{rk}(b,B^*)}(B^*)$ is an atom of
$B^*/\text{id}_{\text{rk}(b,B^*)}(B^*)$
\sn
\itemitem{ $(iii)$ }  if $b' < b \and b' \in \text{id}_{\text{rk}(b,B^*)}
(B^*)$ then $\{x \in {\Cal U}:x \le b'\}$ is finite
\endroster}
\item "{$(k)$}"  if $\lambda' > \mu$ then $\chi,{\Cal A}'$ satisfies:
{\roster
\itemitem{ $(\alpha)$ }  ${\Cal A}' \subseteq [\lambda']^{\aleph_0}$ is a MAD family
of cardinality $\le \chi$ such that: every permutation $\pi$ of $\lambda' \backslash
Z,Z \in [\lambda']^{< \sigma}$ satisfying $A \in {\Cal A}' \Rightarrow |A
\Delta \pi''(A) \backslash Z| < \aleph_0$, has support 
$\{\alpha < \lambda':\pi(\alpha) \ne \alpha\}$ of cardinality $<
\theta$
\sn
\itemitem{$ (\beta)$ }  for some ideal $I^*$ of $B^*$ containing
id$_1(B^*)$ the Boolean algebra $B^*/I^*$ is isomorphic to $\{a
\subseteq \chi:a$ is finite or co-finite$\}$.
\endroster}
\ermn
\ub{Then} we can find $B$ such that:
\mr
\item "{$(\alpha)$}"  $B$ is a superatomic Boolean Algebra
\sn
\item "{$(\beta)$}"  $B$ has $\lambda'$ atoms and $\lambda$ elements
\sn
\item "{$(\gamma)$}"  every automorphism $g$ of $B$ moves $< \theta$
atoms; i.e. \nl
$|\{x \in {\text{ \rm atom\/}}(B):g(x) \ne x\}| < \theta$.
\endroster
\endproclaim
\bigskip

\demo{Proof}  Without loss of generality $B^*$ is a Boolean Algebra of
subsets of $\mu$ with $\{\{\alpha\}:\alpha < \mu\}$ being the atoms of $B^*$.
If $\lambda' = \mu$ let ${\Cal A}'= \emptyset$.

Without loss of generality $B_\alpha$ is a subalgebra of
${\Cal P}(X_\alpha)$ and the set of atoms of $B_\alpha$ is $\{\{x\}:
x \in X_\alpha\}$.  Without loss of generality $\alpha \ne \beta \Rightarrow
X_\alpha \cap X_\beta = \emptyset$ and let $X = \cup\{X_\alpha:\alpha < \mu\}$. 

Let $Y^* \subseteq B^*$ be such that $\{y/I^*:y \in Y^*\}$ is the set
of atoms of $B^*/I^*$ with no repetitions; \wilog \, for each $y \in
Y^*$ for some $\alpha,y/\text{id}_\alpha(B^*)$ is an atom of
$B^*/\text{id}_\alpha(B^*)$ and $(\forall z)[z \le_{B^*} y \rightarrow
z \in \text{ id}_\alpha(B^*) \equiv z \in I^*]$ (just possibly
decrease each $y \in Y^*$).

Let $Y$ be such that $Y \subseteq B^*,\langle
y/\text{ id}_{\text{rk}(y,B^*)}(B^*):y \in Y \rangle$ list with no
repetitions $\{y/\text{ id}_{\text{rk}(y,B^*)}:y/\text{id}_{\text{rk}(y,B^*)}
(B^*)$ an atom of $B^*/\text{id}_{\text{rk}(y,B^*)}(B^*)$ and rk$(y,B^*) 
> 0\}$ and let $D_y$ be the ultrafilter on $B^*$ generated by
$\{y\} \cup \{1-x:x \in B^*,\text{rk}(x,B^*) < \text{ rk}(y,B^*)\}$ for
each $y \in Y$.  Without loss of generality $Y^* \subseteq Y$ also
clearly $y \in Y \Rightarrow \{y' \in Y^*:y'-y \in \text{
id}_{\text{rk}(y',B^*)}(B^*)\}$ is finite so \wilog \, is empty for $y
\in Y \backslash Y^*$ (singleton for $y = Y^*$ of course), note that
$Y^*$ is of cardinality $|{\Cal A}'|$ and 
without loss of generality $|Y \backslash Y^*| = \lambda$.
\nl
First assume
\mr
\item "{$\boxtimes_1$}"  $y \in Y \Rightarrow |y| = \mu$.
\ermn

Let $g$ be a one-to-one function from $\mu$ onto $X$ and for 
$A \in {\Cal A}$ (from clause (b)) let $\{\gamma_{A,k}:k < \omega\}$ 
list $A$ without repetition.  Let $g^*:\mu \rightarrow \mu$ be $g^*(\gamma) =
\text{ Min}\{\alpha < \mu:g(\gamma) \in X_\alpha\}$.  For each 
$A \in {\Cal A}$, choose if possible an infinite $u = u_A \subseteq \omega$ such that
$\langle g^*(\gamma_{A,k}):k \in u \rangle$ is with no 
repetitions and $\langle g^*(\gamma_{A,k}):k \in u \rangle$ converge to 
some $D_y,y = y_A \in Y$.  
Note that the only case $u$ is not well defined, is when the set
$\{g^*(\gamma_{A,k}):k \in u\}$ is finite; we use clause (j) and 
properties of superatomic Boolean Algebras.
As [$y_A$ well define $\Rightarrow |y_A| = \mu$] by $\boxtimes$, clearly we
can find $\langle \alpha[A]:A \in {\Cal A},u_A$ well defined$\rangle$ such
that: $\alpha[A] < \mu,\alpha[A] \in y_A$ and $(\forall z \in Y)[\alpha[A]
\in z \Leftrightarrow z \in D_{y_A}]$ and 
$\alpha[A_1] = \alpha[A_2] \Rightarrow A_1 = A_2$.  
Let for $\alpha < \mu,a_\alpha$ be $\{g(\gamma_{A,k}):k \in u_A\}$ 
if $A \in {\Cal A},\alpha[A] = \alpha$ and $u_A$ 
is well defined, and $\emptyset$ otherwise.
Toward defining our Boolean Algebra let 
$\{x^*_\gamma:\gamma \in [\mu,\lambda')\}$ be
pairwise distinct elements not in $X$.  Let 
${\Cal A}'' = \{\{\mu + i:i \in A\}:A \in {\Cal A}'\}$ so it is a 
maximal almost disjoint family of countable subsets of $[\mu,\lambda')$,
as in clause (k) of the assumption so if $\mu = \lambda'$ then ${\Cal
A}'' = \emptyset = {\Cal A}',\lambda' - \mu = 0,(\lambda'-
\mu)^{\aleph_0} = 0$.

Now we define our Boolean Algebra $\bold B$.  It is the Boolean Algebra of
subsets of $X^* = \dbcu_{\alpha < \mu} X_\alpha \cup \{x^*_\gamma:\gamma \in
[\mu,\lambda')\}$ generated by the following:
\mr
\widestnumber\item{$(iii)$}
\item "{$(i)$}"   the sets $\{a \in B_\alpha:|a| < \theta\} \cup
\{a \cup a_\alpha:a \in B_\alpha,|a| \ge \theta\}$ when $\alpha < \mu$
\sn
\item "{$(ii)$}"  $\{x^*_\gamma\}$ for $\gamma \in [\mu,\lambda')$
\sn
\item "{$(iii)$}"  the sets $c_y$ (for $y \in Y$) where \nl
$$
\align
c_y = \bigl\{ x \in X:&\text{for some } \alpha < \mu \text{ we have } 
x \in X_\alpha \and \{\alpha\} \le_{B^*} y\} \cup \\
  &\{x^*_\gamma:\gamma \in [\mu,\lambda') \text{ and } y \in Y^*
\text{ and } \gamma \in d_y \bigr\}.
\endalign
$$
\ermn
Clearly
\mr
\item "{$\otimes_1$}"  $\bold B$ is a subalgebra of ${\Cal P}(X^*)$, including
all the singletons hence is atomic; has $\lambda'$ atoms and $\lambda$
elements
\sn
\item "{$\otimes_2$}"  for $\alpha < \mu$, we have $a \in B_\alpha
\and |a| < \theta \Rightarrow a \in \bold B \and \bold B \restriction a =
B_\alpha \restriction a$ but $a \in B_\alpha \Rightarrow B_\alpha
\restriction a$ is superatomic so $\{a \in B_\alpha:|a| < \theta\}
\subseteq \text{ id}_\infty(\bold B)$. \nl
[Why?  For the first implication we should check that every one of the
generators of $\bold B$ lited in $(i), (ii), (iii)$ above satisfies:
its intersection with $a$ belong to $B_\alpha \restriction a$.  The
rest follows.]    
\sn
\item "{$\otimes_3$}"  for $\alpha < \mu,\{a \in \bold B:a
\subseteq X_\alpha \cup a_\alpha:|a| < \theta\}$ satisfies
{\roster
\itemitem{ $(i)$ }  it is equal to $\{a \cup b:a \in B_\alpha \and |a| <
\theta \text{ and } b \subseteq a_\alpha \text{ is finite}\}$
\sn
\itemitem{ $(ii)$ }  it is a maximal ideal of $\bold B \restriction (X_\alpha
\cup a_\alpha)$
\sn
\itemitem{ $(iii)$ }  it is included in id$_\infty(\bold B)$. \nl
[Why?  Just think.]
\endroster}
\item "{$\otimes_4$}"  $\alpha < \mu \Rightarrow X_\alpha \cup a_\alpha 
\in \text{ id}_\infty(\bold B)$ \nl
[Why?  First $X_\alpha \cup a_\alpha \in \bold B$ by clause (i) above,
second if $X_\alpha \cup a_\alpha \notin \text{ id}_\infty(\bold B)$
then by $\otimes_3$ above $(X_\alpha \cup a_\alpha)$ is an atom of
$\bold B/\text{id}_\alpha(B)$ for $\alpha$ large enough, hence
$X_\alpha \cup a_\alpha$ belong to id$_{\alpha +1}(\bold B)$,
contradiction
\sn
\item "{$\otimes_5$}"  for $\alpha < \mu,\bold B \restriction (X_\alpha \cup
a_\alpha) \cong B_\alpha$ hence if $\alpha < \beta < \omega$ then for no
$c_\alpha$ such that $c_\alpha \in B_\alpha,
c_\alpha \le X_\alpha \cup a_\alpha,|c_\alpha| < 
\theta$ and $c_\beta \in B_\beta,c_\beta \le X_\beta \cup a_\beta,
|c_\beta| < \theta$ do
we have $\bold B \restriction (X_\alpha \cup a_\alpha \backslash c_\alpha)
\cong \bold B \restriction (X_\beta \cup a_\beta \backslash c_\beta)$. \nl
[Why?  By clauses (f) + (e) of the assumption, the ``hence" follows by
clause (g) of the assumption.]
\ermn
Let $I_1$ be the ideal $[X^*]^{< \theta} \cap \bold B$ of $\bold B$.
So clearly
\mr
\item "{$\otimes_6$}"  $I_1 \subseteq \text{ id}_\infty(\bold B)$.
\ermn
We shall prove that
\mr
\item "{$\otimes_7$}"  $\bold B/I_1$ is isomorphic to a homomorphic 
image of $B^*$.
\ermn
Toward proving $\otimes_7$ let $S = \{x^*_\gamma:\gamma \in [\mu,\lambda')\}$
and define a function $h$ as follows: its domain is $\{c_y:y \in Y\} \cup
\{X_\alpha \cup a_\alpha:\alpha < \mu\}$ and $h(c_y) = y,h(X_\alpha \cup
a_\alpha) = \{\alpha\}$.  Now
\mr
\item "{$(*)_0$}"  $(X_\alpha \cup a_\alpha)/I_0$ is an atom of
$\bold B/I_1$ \nl
[why?  by $\otimes_3$.]
\sn
\item "{$(*)_1$}"  $\{b/I_1:b \in \text{ Dom}(h)\}$ is a subset of
$\bold B/I_1$ which generates it \nl
[why?  see the definitions of $\bold B$ and of $I_1$.]  
\sn
\item "{$(*)_2$}"  if $n_1 \le n < \omega,m_1 \le m < \omega,y_0,\dotsc,
y_{n-1} \in Y$ is with no repetitions, $\alpha_0,\dotsc,\alpha_{m-1} < \mu$
is with no repetitions, \ub{then}: \nl
in $\bold B,\tau_1 =: \dbca_{\ell < n_1} c_{y_\ell} \cup \dbca_{\ell < m_1}
(X_{\alpha_\ell} \cup a_{\alpha_\ell}) - \dbcu^{n-1}_{\ell =n_1} \,
c_{y_\ell} \cup \dbcu^{m-1}_{\ell = m_1} (X_{\alpha_\ell} \cup
a_{\alpha_\ell})$ belongs to $I_1$ \ub{iff} \nl
in $B^*,\tau_2 =: 
\dbca_{\ell < n_1} y_\ell \cup \dbca_{\ell < m_1}\{\alpha_\ell\}
- \dbcu^{n-1}_{\ell = n_1} \, y_\ell \cup \dbcu^{m-1}_{\ell = m_1} 
\{\alpha_\ell\}$ is empty. \nl
[Why?  First, assume that the second statement holds then trivially
$\tau'_1 =: \dbca_{\ell < n_1} (c_{y_\ell} \backslash S) \cup 
\dbca_{\ell < m_1} X_{\alpha_\ell} - \dbcu^{n-1}_{\ell =n_1}
(c_{y_\ell} \backslash S) \cup \dbcu^{m-1}_{\ell = m_1} X_{\alpha_\ell} =
\bigcup \{X_\beta:B^* \models \{\alpha\} \le \tau_2\} = \emptyset$ but
$\tau'_1 \Delta \tau_1 \subseteq S \cup \dbcu_{\ell < m} 
a_{\alpha_\ell}$ but $a_{\alpha_\ell} \in I_0 \subseteq [X^*]^{< \theta}$ 
and $\tau'_1 = \emptyset$, so $\tau_1 \subseteq S$ mod $[X^*]^{< \theta}$.
\nl
Now if $\tau_1 \cap S \notin J_0$ then $\tau_1 \cap S$ is infinite,
clearly $\lambda' > \mu$, so as $\{d_z:z \in Y^*\}$ is a MAD family of
subsets of $\lambda' \backslash \mu$, necessarily for some $z \in Y^*$
we have $\tau_1 \cap S \cap d_z$ is infinite.  As $\tau_1 \cap S \cap
d_z \subseteq c_{y_\ell}$ for $\ell < n_1$, necessarily $y_\ell = z$,
hence $y_0 = z,n_2=1$.  Similarly $\ell \in [n_1,n_2) \Rightarrow
y_\ell \ne z$ hence $\ell \in [n_1,n) \Rightarrow y_\ell \cap y_0 =
y_\ell \cap z \in \text{ id}_{\text{rk}(z,B^*)}(B^*) \Rightarrow |d_z
\cap c_{y_\ell}| < \aleph_0$.  Hence
clearly $\ell \in [n_1,n) \Rightarrow y_\ell \notin D_0$ but
$y_0 \in D_z$ and $\alpha < \mu \Rightarrow \{\alpha\} \notin D_z$
(as $y \in Y!$) hence $B^* \models ``\tau_2 > 0"$, contradiction,  so
necessarily $\tau_1 \cap S$ is finite hence $\in I_1$. \nl
Second, if the second statement fails, then for some $\beta < \mu,B^*
\models \{\beta\} \le \tau_2$, but then $X_\beta \subseteq \tau'_1$ and as
above $\tau'_1 \supseteq \tau_1 \backslash S \supseteq X_\beta$ mod
$I_1$ but $S \cap X_\beta = \emptyset$, so $X_\beta \subseteq \tau_1$
mod $I_1$; now $X_\beta \notin I_0$ (as $|X_\beta| \ge \theta$ by clause (e)
of the assumptions) hence $\tau_1 \notin I_1$.  So we have proved $(*)_2$.]
\ermn
Now by $(*)$, $\otimes_7$ follows, in fact $h$ induces an isomorphism
$\hat h$ from $\bold B/I_0$ onto $B^*$.  But $B^*$ is superatomic and $I_0
\subseteq \text{ id}_\infty(\bold B)$ by $\otimes_6$ hence
\mr
\item "{$\otimes_8$}"  $\bold B$ is superatomic.
\ermn
For the rest of the proof let $f \in \text{ AUT}(\bold B)$ and toward
contradiction we assume sup$(f) = \{x \in \text{ atom}(\bold B):f(x) \ne
x\}$ has cardinality $\ge \theta$. \nl
Recall that $I_1 = \{a \in \bold B:|a| < \theta\}$ so necessarily $f$ maps $I_1$ 
onto itself.
Note that $\{X_\alpha \cup a_\alpha/I_1:\alpha < \mu\}/I_1$ list the atoms
of $\bold B/I_1$.  Assume $f(X_\alpha \cup a_\alpha)/I_1 = (X_\beta \cup
a_\beta)/I_1,\alpha \ne \beta$; let $c_1 = (X_\alpha \cup a_\alpha) -
f^{-1}(X_\beta \cup a_\beta)$ and $c_2 = (X_\beta \cup a_\beta) - f(X_\alpha
\cup a_\alpha)$, so both being the difference of two members of $\bold B$ are
in $\bold B$ and $c_1 \le X_\alpha \cup a_\alpha,c_2 \le X_\beta \cup
a_\beta$ hence $c_1 \in B_\alpha,c_2 \in B_\alpha$.  Clearly $f
\restriction (\bold B \restriction 
(X_\alpha \cup a_\alpha - c_1))$ is an isomorphism from $\bold B \restriction
(X_\alpha \cup a_\alpha - c_1)$ onto $\bold B \restriction (X_\beta \cup
a_\beta - c_2)$, contradicting $\otimes_5$.  Hence the automorphism $f$
induced on $B^*/I_1$ maps each atom to itself hence 
is the identity.  Also for $\alpha < \mu$ we have
$(X_\alpha \cup a_\alpha) \Delta f(X_\alpha \cup a_\alpha) \in I_1$, that 
is, has cardinality $< \theta$.  So for each $\alpha < \mu$, 
letting $c^1_\alpha =: (X_\alpha \cup a_\alpha) - f^{-1}(X_\alpha \cup
a_\alpha) \in I_\alpha$ and $c^2_\alpha =: (X_\alpha \cup a_\alpha) =
f(X_\alpha \cup a_\alpha) \in I_\alpha,f \restriction (B_\alpha
\restriction (1 - c^1_\alpha))$ is an isomorphism from $B_\alpha
\restriction (1 - c^2_\alpha)$ onto $B_\alpha \restriction (1 - c_2)$
hence
\mr
\item "{$\boxtimes_2$}"  $Z_\alpha = \{x:x$ an atom of $B_\alpha,x
\le_{B_\alpha} c^1_\alpha \vee f(x) \ne x\}$ has cardinality $<
\theta$
\ermn
by clause (d) of the assumptions on $B_\alpha$.
Let $v =: \{\alpha < \mu:\text{for some } x \in X_\alpha$ we have
$f(\{x\}) \ne \{x\}\}$.  Assume for the time being
\mr
\item "{$\boxtimes_3$}"  $v$ has cardinality $\ge \text{ cf}(\theta)$.
\ermn
For $\alpha \in v$ choose $x_\alpha \in X_\alpha$ such that
$f(x_\alpha) \ne x_\alpha$ and shrinking $v$ without loss of generality 
$\alpha,\beta \in v \Rightarrow x_\alpha \ne f(x_\beta)$.  
Let $g:v \rightarrow \mu +1$ be such that $f(x_\alpha) \subseteq
X_{g(\alpha)}$ where we stipulate $X_\mu = S$.  Applying the above to
$f^{-1}$ \wilog \, either $g$ is one-to-one into $\mu$ or $g$ is
contantly $\mu$.  Hence by clause (b) of the assumption
for some $A \in {\Cal A}$ we have $(\forall \gamma \in A)[\{\gamma\}
\in \{x_\alpha:\alpha \in v\}$.  So 
$\alpha[A] < \mu$ is well defined and an easy contradiction.
We can conclude that $\neg \boxtimes_3$ hence
$Z =: \{x \in X:f(\{x\}) \ne \{x\}\}$
has cardinality $< \text{ cf}(\theta)$ hence $|\{x \in X:f(x) \ne x\}| <
\theta$.  If $\mu = \lambda'$ we are done so assume $\mu < \lambda'$.
\sn
Now $S = \{x^*_\gamma:\gamma \in [\mu,\lambda')\} = X^* \backslash X
\subseteq X^*$ satisfies:
\mr
\item "{$\otimes_9(\alpha)$}"  $(\forall b \in \bold B)(b \cap S$ infinite $\Rightarrow 1
\le \text{ rk}(b/I_1,\bold B/I_1))$ and
\sn
 \item "{${{}}(\beta)$}"   if $S'$ satisfies the property in clause
$(\alpha)$,  then $|S' \backslash S| < \sigma$ \nl
[Why?  Clause $(\alpha)$ is proved by inspecting the definition of
$\bold B$.  As for clause $(\beta)$, if
$|S' \backslash S| \ge \sigma$ as $S' \backslash S \subseteq X$
clearly then there is $A \in {\Cal A}$ such that $\{g(i):i \in A\} \subseteq
S' \backslash S$.  First if $\alpha =: \alpha[A]$ is well defined then $X_\alpha \cup
a_\alpha \in \bold B$, rk$((X_\alpha \cup a_\alpha)/I_1,\bold B/I_1) = 0 < 1$
but $(X_\alpha \cup a_\alpha) \cap S' \supseteq a_\alpha$ is infinite;
contradiction.  Second if $\alpha[A]$ is not well defined then for some $\alpha
< \mu$ we have $\{g(i):i \in A\} \cap X_\alpha$ is infinite and we get a
similar contradiction.]
\ermn
Hence $S^*_f =: \{x^*_\gamma:\gamma \in [\mu,\lambda')$ and 
$f^n(\{x^*_\gamma\}) \subseteq X$ or $f^{-n}(\{x^*_\gamma\}) \subseteq
X$ for some $n < \omega\}$ 
has cardinality $< \sigma$ (this also follows from the
previous paragraph recalling $\sigma = \text{ cf}(\sigma) > \aleph_0$).

Also for $y \in Y^*$ letting $\gamma = \text{ rk}(y,B^*)$ we have
$c_y \Delta f(c_y) \in I_1$, (just recall that the automorphism $f$
induced on $\bold B/I_1$ is the identity, and recall that $[d \subseteq S
\and d \in I_1 \Rightarrow d$ is finite], hence the symmetric difference of
$\{\{x^*_\gamma\}:\gamma \in d_y\} \backslash S^*_f,
\{f\{x^*_\gamma\}:\gamma \in d_y\} \backslash S^*_f$ is finite. \nl
As $\{d_y:y \in Y^*\}$ is MAD as inclause the set 
$\{\gamma \in [\mu,\lambda'):f(\{x^*_\gamma\}) \ne \{x^*_\gamma\}\}$ 
is of cardinality $< \theta$; so seemingly we are done.  
\enddemo
\bn
Not exactly: we have assumed $\boxtimes$, i.e. 
$y \in Y \Rightarrow |y| \ge \mu$.

To eliminate this we make some minor changes.  First \wilog \, $B^*$ is a
Boolean Algebra of subsets of $\{\alpha:\alpha < \mu \text{ even}\}$ with 
the singletons being its atoms.  Second, for $A \in {\Cal A}$, if possible
we choose $u = u_A$ as follows:
\mr
\item "{$(a)$}"  either $(\alpha)$ or $(\beta)$ where
{\roster
\itemitem{ $(\alpha)$ }  $g^*(\gamma_{A,k})$ is odd for every $k \in u$
\sn
\itemitem{ $(\beta)$ }  $g^*(\gamma_{A,k})$ is even for every $k \in u$
\endroster}
\item "{$(b)$}"  $\langle g^*(\gamma_{A,k}):k \in u \rangle$ is with no
repetitions
\sn
\item "{$(c)$}"  if case $(\beta)$ occurs in $A$, then there is a unique
$y = y_A \in Y$ such that $\langle \{g^*(\gamma_{A,k})\}:k \in u \rangle$
converge to $D_{y_A}$.
\ermn
Note
\mr
\item "{$(*)$}"  if $u_A$ is not well defined then for some finite $w 
\subseteq \mu$ we have \nl
$\{g(\gamma_{A,k}):k < \omega\} \subseteq
\dbcu_{\alpha \in \omega} X_\alpha$.
\ermn
Now we choose $\langle \alpha[A]:A \in {\Cal A},u_A$ will define $\rangle$
such that:
\mr
\item "{$(**)$}"  $\langle \alpha[A]:A \in {\Cal A},u_A$ well defined
$\rangle$ list with no repetitions the odd ordinals $< \mu$.  
\ermn
We define $a_\alpha$, etc. as before.
\sn
Lastly, defining $c_y$, we add $\{x$:\, for some $A \in {\Cal A},\alpha[A]$ is
well defined, $x \in X_{\alpha[A]}$ and $B^* \models y_A \le y\}$. \nl
Note that we just replace $B^*$ by $B^{**} \subseteq {\Cal P}(\mu)$
where \wilog \, $B^* \subseteq {\Cal P}(\mu)$, atom$(B^*) =
\{\{\alpha\}:\alpha < \mu\}$, let $f:\mu \rightarrow \mu$ be onto such
that $(\forall \alpha)(\exists^\mu \beta)[f(\beta) = \alpha]$ and we
let $B^{**}$ be the subalgebra of $\bold P(\mu)$ generated by
$\{\{\alpha\}:\alpha < \mu\} \cup \{\{\beta < \mu:f(\beta) \in y\}:y
\in B^*\}$.
\hfill$\square_{\scite{gd.2}}$
\bn
\ub{Discussion}:  

Why do we use MAD families ${\Cal A} \subseteq [\mu]^{\aleph_0}$ and not
$\subseteq [\mu]^{\aleph_1}$?  If we use the latter, we have to take more
care about superatomicity as the intersections of such members may otherwise
contradict superatomicity.
\newpage

\head {\S3 Sufficient conditions for the existence of $\langle B_i:i < \mu 
\rangle$} \endhead  \resetall 
\bigskip

Here we shall show that the assumptions of \scite{gd.1} are reasonable.  Now
in \scite{ge.0a} we shall reduce the clause (k) of \scite{gd.1} to Pr$(\lambda',
\theta)$ where Pr formalizes clause (b) there.  In \scite{ge.1}, \scite{ge.3}
we give sufficient conditions for Pr$(\mu,\sigma)$.  In fact, it is clear
that (high enough) it is not easy to fail it.  In \scite{ge.7} we give a
sufficient condition for a strong version of clauses (e) - (f) of
\scite{gd.1} (and earlier deal with the conditions appearing in it).
So for not having the assumptions of \scite{gd.1} it has large
consistency strength.
\bigskip

\definition{\stag{ge.0} Definition}  1) Pr$(\chi,\mu,\sigma)$ means that for
some ${\Cal A}$ we have:
\mr
\item "{$(a)$}" ${\Cal A} \subseteq [\mu]^{\aleph_0}$
\sn 
\item "{$(b)$}"  ${\Cal A}$ is almost disjoint, i.e. $A \ne B \in {\Cal A}
\Rightarrow |A \cap B| < \aleph_0$
\sn
\item "{$(c)$}"  $|{\Cal A}| = \chi$
\sn
\item "{$(d)$}"  $(\forall A \in [\mu]^\sigma)(\exists A \in {\Cal A})
[A \subseteq^* A]$.
\ermn
2) If we omit $\chi$ we mean ``some $\chi$". \nl
3) We call ${\Cal A} \subseteq [\lambda]^{\aleph_0}$ saturated if for every
$A \in [\lambda]^{\aleph_0}$ not almost contained in a finite union of
members of ${\Cal A}$, almost contains a member of ${\Cal A}$.
\enddefinition
\bn
\ub{\stag{ge.0a} Fact}:  1) Clause (b) of the assumption of \scite{gd.2} is
equivalent to Pr$(\mu,\mu,\sigma)$. \nl
2) Clause (k)$(\alpha)$ of the 
assumption of \scite{gd.2} follows from Pr$(\chi,\lambda',\sigma) \and
\chi \ge 2^{\aleph_0}$. \nl
3) If ${\Cal A} \subseteq [\mu]^{\aleph_0}$ is almost disjoint and is saturated
\ub{then} Pr$(|{\Cal A}|,\mu,\aleph_1)$. \nl
4) If $\mu = \mu^{\aleph_0} \ge \sigma$ then Pr$(\mu,\sigma) \equiv
\text{ Pr}(\mu,\mu,\sigma)$ and $\chi \ne \mu \Rightarrow \neg \text{ Pr}
(\chi,\mu,\sigma)$. \nl
5) For any $\lambda' \ge \aleph_0$ \ub{there} is a MAD family 
${\Cal A} \subseteq [\lambda']^{\aleph_0}$ of cardinality 
$[\lambda']^{\aleph_0}$ satisfying clause (k)$(\alpha)$ of \scite{gd.2}.
\bigskip

\demo{Proof}  1) Read the two statements. \nl
2) Let ${\Cal A} \subseteq [\lambda']^{\aleph_0}$ exemplify Pr$(\lambda',
\sigma)$.  For each $A \in {\Cal A}$ we can find $\langle B_{A,\zeta}:
\zeta < 2^{\aleph_0} \rangle$ such that:
\mr
\widestnumber\item{$(iii)$}
\item "{$(i)$}"  $B_{A,\zeta} \in [A]^{\aleph_0}$
\sn
\item "{$(ii)$}"  $\zeta \ne \varepsilon \Rightarrow B_{A,\zeta} \cap
B_{A,\varepsilon}$ is finite
\sn
\item "{$(iii)$}"  if $\pi$ is a partial one-to-one function from $A$ to
$A$ such that $x \in \text{ Dom}(\pi) \rightarrow x \ne \pi(x)$
\ub{then}
for some $\zeta < 2^{\aleph_0}$ we have $\alpha \in B_{A,\zeta}
\Rightarrow \alpha \notin \text{ Dom}(\pi) \vee \pi(\alpha) \notin
B_{A,\zeta}$. 
\ermn
Why? First find $\langle B'_{A,\zeta}:\zeta < 2^{\aleph_0} \rangle$ satisfying
(i), (ii), let $\langle \pi_\zeta:\zeta < 2^{\aleph_0} \rangle$ list the
$\pi$'s from (iii) and chose $B_{A,\zeta} \in [B'_{A,\zeta}]^{\aleph_0}$ to
satisfy clause (iii) for $\pi_\zeta$.
\nl
Lastly, ${\Cal A}' = \{B_{A,\zeta}:A \in {\Cal A}$ and $\zeta < 2^{\aleph_0}
\}$ it is as required.  \nl
3), 4)  Easy. \nl
5) Starting with AD ${\Cal A}_0 \subseteq [\lambda']^{\aleph_0}$ of
cardinality $9\lambda')^{\aleph_0}$, extend it to a MAD one ${\Cal
A}_1$ and then apply the proof of part (2).
\hfill$\square_{\scite{ge.0a}}$
\enddemo
\bigskip

\proclaim{\stag{ge.1} Claim}  1) Assume
\mr
\item "{$(a)$}"  $\kappa_n < \kappa_{n+1} < \kappa < \mu_n < \mu_{n+1} <
\mu$ for $n < \omega$
\sn
\item "{$(b)$}"  $\kappa = \sum \kappa_n,\mu = \sum \mu_n$ and {\rm max pcf}
$\{\kappa_n:n < \mu\} > \mu$
\sn
\item "{$(c)$}"  $\kappa$ strong limit and $2^\kappa \ge \mu^+$
\sn
\item "{$(d)$}"  $\langle \mu_n:n < \omega \rangle$ satisfies the
requirements from \cite[\S1]{Sh:513} or at least the conclusion.
\ermn
\ub{Then} for every $\lambda \ge \kappa$ we can find $\{\bar A_\alpha:\alpha
< \alpha^*\}$ such that
\mr
\item "{$(\alpha)$}"   each $\bar A_\alpha$ has the form $\langle
A_{\alpha,n}:n < \omega \rangle$, it  belongs to $\dsize \prod_{n < \omega}
[\lambda]^{\kappa_n}$ and for each $\alpha$ we have $\langle A_{\alpha,n}:
n < \omega \rangle$ pairwise disjoint
\sn
\item "{$(\beta)$}"  if $\alpha \ne \beta$, then $\bar A_\alpha,A_\alpha$ are
almost disjoint which means $f \in \dsize \prod_{n < \omega}
A_{\alpha,n} \and f' \in \dsize \prod_{n < \omega} A_{\beta,n} \Rightarrow
|{\text{\rm Rang\/}}(f) \cap {\text{\rm Rang\/}}(f')| < \aleph_0$
\sn
\item "{$(\gamma)$}"  if $\bar A \in \dsize \prod_{n < \omega} [\lambda]
^{\kappa_n}$ \ub{then} for some $\alpha < \alpha^*$ and
one to one function $h_1,h_2 \in 
{}^\omega \omega$ we have $\kappa = {\text{\rm lim\/}} \langle|A_{h_1(n)} \cap
A_{\alpha,h_2(n)}|:n < \omega \rangle$.
\ermn
2) We can conclude in (1) that: there is ${\Cal A} \subseteq [\lambda]
^{\aleph_0}$, an almost disjoint family such that $(\forall B \in
[\lambda]^\kappa)(\exists A \in {\Cal A})(A \subseteq B)$. \nl
3) If in part (1) instead (a)-(d) we just assume $\kappa$ strong limit
of cofinality $\aleph_0$, and SCH + $(\forall \lambda >
\kappa)[\text{cf}(\lambda) = \aleph_0 \rightarrow
\diamondsuit_{\lambda^+}$, \ub{then} the conclusion of (1) holds.
\endproclaim 
\bigskip

\demo{Proof}  By \cite{Sh:460}, \cite[\S3]{Sh:668} (even more).  
\enddemo
\bigskip

\remark{\stag{ge.2} Remark}  Are the hypotheses of \scite{ge.1}(1)
reasonable? \nl
1)  If for some strong limit $\kappa$ of cofinality $\aleph_0,2^\kappa >
\kappa^{+ \omega}$, then we can let $\mu_n = \kappa^{+1+n}$ and $\langle 
\kappa_n:n < \omega \rangle$ as in clause (a), (b), (c) there exists (by
\cite[Ch.IX,\S5]{Sh:g}, and it is hard not to satisfy clause (d)
(see \cite{Sh:513}). \nl
2) Clause (c), i.e. $\kappa$ strong limit, is needed just to start the
induction.  If $\kappa = \aleph_\omega \le 2^{\aleph_0}$ we have a similar
theorem. \nl
3) If $\aleph_\omega \le 2^{\aleph_0}$ and $\aleph_\omega$ is as required
in \cite[\S1]{Sh:513} then we have a parallel theorem.
\endremark
\bn
We quote \cite{GJSh:399} in \scite{ge.3}(1) and (2) is immediate
starting the induction with the known Pr$(\lambda,\aleph_1)$ for
$\aleph_0 < \lambda \le 2^{\aleph_0}$.
\proclaim{\stag{ge.3} Claim}  1) Assume CH + SCH + $(\forall \mu)
({\text{\rm cf\/}}(\mu) = \aleph_0 < \mu \rightarrow \square_{\mu^+})$.
\ub{Then} there is a saturated MAD family ${\Cal A}_\lambda \subseteq
[\lambda]^{\aleph_0}$ for every uncountable $\lambda$. \nl
2) If SCH + $(\forall \mu)[\text{cf}(\mu) = \aleph_0 \rightarrow
\square_{\mu^+})$, \ub{then} Pr$(\lambda,\aleph_1)$ for every $\lambda > \aleph_0$.
\endproclaim
\bigskip

\definition{\stag{ge.4} Definition}   Let $\mu \ge \theta$. \nl
1) Let ${\Cal S}_\theta$ be the class of $\bar a = \langle a_n:n < \omega
\rangle$ such that $|a_n| \le \theta$, cf$(\theta) = \aleph_0 \Rightarrow
|a_n| < \theta$ and $\theta = \text{ lim sup}_n|a_n|$.  Let 
${\Cal S}_{\theta,\mu} =: \{\bar a:\bar a =
\langle a_n:n < \omega \rangle,a_n \in [\mu]^{\le \theta},a_n \subseteq
a_{n+1}$ and $\theta = \text{ lim sup}_{n < \omega}|a_n|\}$. \nl
3) For $\bar a \in {\Cal S}_\theta$ let set$(\bar a) =
\{w:|w| = \aleph_0$ and $w \subseteq \dbcu_{n < \omega} a_n$ and $n < \omega
\Rightarrow |w \cap a_n \backslash \dbcu_{\ell < n} a_\ell| < \aleph_0\}$.
\nl
4) For $\bar a,\bar b \in {\Cal S}_\theta$ let $\bar a \le^* \bar b$ mean
set$(\bar a) \supseteq \text{ set}(\bar b)$. \nl
5) We say 
$\bar a,\bar b \in {\Cal S}_\theta$ are compatible if 
$(\exists \bar c \in {\Cal S}_\theta)[\bar a \le^* \bar c \and \bar b
\le^* \bar c \and \dbcu_n c_n \subseteq \dbcu_n a_n \cap \dbcu_n b_n]$.
\enddefinition
\bigskip

\definition{\stag{ge.5} Definition}  Let $\boxtimes_{\theta,\mu}$ be
\roster
\item "{$\boxtimes_{\theta,\mu}$}"  there is ${\Cal S}^* \subseteq
{\Cal S}_{\theta,\mu}$ such that:
{\roster
\itemitem{ $(a)$ }  for every 
$\bar a \in {\Cal S}$ there is $\bar b \in {\Cal S}^*$ compatible with
$\bar a$
\sn
\itemitem{ $(b)$ }  if $\bar a \ne \bar b \in {\Cal S}^*$ then the 
set$(\bar a) \cap$ set$(\bar b) = \emptyset$.
\endroster}
\endroster
\enddefinition
\bigskip

\proclaim{\stag{ge.6} Claim}  If $\theta$ is strong limit, $\theta >
\text{ cf}(\theta) = \aleph_0$ and $\mu \in (\theta,(2^\theta)^{+ \omega})$
satisfies $\otimes_{\theta,\mu}$ below \ub{then}
$\boxtimes_{\theta,\mu}$ from
\scite{ge.5} holds where:
\mr
\item "{$\otimes_{\theta,\mu}$}"  $\mu = 2^\theta$, 
{\rm pp}$_{J^{\text{\rm bd\/}}_\omega}(\theta) = 2^\theta$.
\endroster
\endproclaim
\bigskip

\demo{Proof}  Straight.  First assume $\mu \le 2^\theta$ as ${\Cal
S}_{\theta,\mu} = \mu^\theta = 2^\theta$, we can find $\langle \bar
a^\alpha:\alpha < 2^\theta \rangle$ listing ${\Cal S}_{\theta,\mu}$.
Now we choose $\gamma_0(\alpha),\bar b^\alpha$ by induction on $\alpha
< 2^\theta$ such that
\mr
\item "{$(a)$}"  $\bar b^\alpha \in {\Cal S}_{\theta,\mu}$
\sn
\item "{$(b)$}"  $\beta < \alpha \Rightarrow \text{ set}(\bar a) \cap
\text{ set}(\bar b) = \emptyset$
\sn
\item "{$(c)$}"  $\bar a^{\gamma(\alpha)},\bar b^\alpha$ are
compatible
\sn
\item "{$(d)$}"   $\gamma(\alpha) = \text{ Min}\{\gamma:\bar a^\gamma$
incompatible with $\bar b^\beta$ for every $\beta < \alpha\}$.
\ermn
Arriving to $\alpha$ choose $\gamma(\alpha)$ by clause (d), we can
find $\kappa_n = \text{ cf}(\kappa_n) < \theta$ such that $\dsize
\prod_{n < \omega} \kappa_n/J^{\text{bd}}_\omega$ has true cofinality
$> |\alpha|$.  Let $h:\omega \rightarrow \omega$ be increasing such
that $|b^{\gamma(\alpha)}_{k(n)}| \ge \kappa_n$ choose
$\gamma_{\alpha,\zeta,n} \in b^{\gamma(\alpha)}_{h^n}$ for $\zeta <
\kappa_n$ increasing with $\zeta$.  For each $\beta < \alpha$ define
$g_{\beta,\alpha} \in \dsize \prod_n(\kappa^{+1}_n)$ by
$g_{\beta,\alpha}(n) = \text{ sup}\{\zeta:\gamma_{\alpha,\zeta,n} \in
\dbcu_{m < \omega} a^\beta_m\}$.  Easily $g_{\beta,\alpha}
<_{J^{\text{bd}_\omega}} \langle \kappa_n:n < \omega \rangle$ hence
there is $\langle \zeta_n:n < \omega \rangle$ such that $(\forall
\beta < \alpha)[g_{\beta,\alpha} <_{J^{\text{bd}_\omega}} \langle
\zeta_n:n < \omega \rangle]$ and let $b^\alpha_n =
\{\gamma_{\alpha,\zeta,n}:n < \omega$ and $\zeta \in
[\zeta_n,\kappa_n)\}$.  Easy to carry and give the conclusion.  For
$\mu = (2^\theta)^{+n}$, use induction on $n$ (as in \cite[\S3]{Sh:668}
or \cite[p.xx]{EH:1}. 
\enddemo
\bigskip

\proclaim{\stag{ge.7} Claim}  1) Assume $\theta$ is strong limit, 
$\aleph_0 = \text{ cf}(\theta) < \theta$.  Assume further 
$\theta \le \kappa \le 
2^{2^\theta},\mu = 2^\kappa$ and $\boxtimes_{\theta,\kappa}$ (from
\scite{ge.5}) holds and $\mu = \mu^{\aleph_0}$.
\ub{Then} some $\bar B = \langle B_\alpha:\alpha < \mu \rangle$ 
satisfies clauses (b) - (g) of \scite{gd.2}; in fact $B_\alpha$ is a subalgebra of ${\Cal P}
(2^\theta)$ with 2 levels and id$_{< \infty}(B_\alpha)$ is included 
is $\{a \subseteq 2^\theta$: a countable or co-countable$\}$. \nl
2) As above except that $2^\theta = \theta^{\aleph_0},\theta > \text{ cf}
(\theta) = \aleph_0$.
\endproclaim
\bigskip

\demo{Proof}  1) Let for simplicity $\theta = \sum \theta_n,\theta_n <
\theta_{n+1} < \theta$.
\bn
\ub{Fact}:  Letting $\bar a^* = \langle \theta_n:n < \omega \rangle$, i.e.
$a^*_n = \theta_n$ we can find $\bar t^{\bar a} = \langle
t_{\ell,\alpha}:\ell < 3,\alpha < 2^\theta \rangle$ such that:
\mr
\widestnumber\item{$(iii)$}
\item "{$(i)$}" $t_{\ell,\alpha} \in$ set$(\bar a^*)$ has order type $\omega$
\sn
\item "{$(ii)$}"  for some one to one onto $\pi:2^\theta \times 2^\theta
\rightarrow 2^\theta$ we write $t_{2,\alpha,\beta}$ for $t_{2,\pi(\alpha,
\beta)}$
\sn
\item "{$(iii)$}"  if $(\ell_1,\alpha_1) \ne (\ell_2,\alpha_2)$ \ub{then}
$t_{\ell_1,\alpha_1} \cap t_{\ell_2,\alpha_2}$ is finite
\sn
\item "{$(iv)$}"   if $\bar a \in {\Cal S}_{\theta,\kappa}$ and
$\dbcu_{n < \omega} a_n \subseteq \theta$ \ub{then} 
for some $\alpha < 2^\theta$ we have
$\beta < 2^\theta \Rightarrow t_{2,\alpha,\beta} \in \text{ set}(\bar a)$ 
\sn
\item "{$(v)$}"  if $\bar a,\bar b \in {\Cal S}_{\theta,\kappa}$ and
$\dbcu_{n < \omega} a_n \cup \dbcu_{n < \omega} b_n \subseteq \theta$ and 
set$(\bar a) \cap$ set$(\bar b) = \emptyset$ and 
$h:\dbcu_{n < \omega} a_n \rightarrow
\dbcu_{n < \omega} b_n$ is one to one and maps $a_n$ onto $b_n$ \ub{then}
for some $\alpha,t_{0,\alpha} \in$ set$(\bar a),t_{1,\alpha} \in$ 
set$(\bar b)$ and $h$ maps $t_{0,\alpha}$ into a co-infinite subset of
$t_{1,\alpha}$.
\endroster
\enddemo
\bigskip

\demo{Proof of the fact}  Straight.
\enddemo
\bn
\ub{Construction}:  Let ${\Cal S}^* = \{\bar a^\gamma:\gamma < \kappa\}$
exemplify $\boxtimes_{\theta,\kappa}$, \wilog \, $\bar a \in {\Cal S}^*
\Rightarrow \dsize \bigwedge_{n < \omega} |a_n| = \theta_n$; (because for
every $\bar a \in {\Cal S}_\theta$ there is $\bar a' \in S_\theta,|a'_n| =
\theta_n$ and set$(\bar a') = \text{ set}(\bar a))$.  Let
$\{X_\gamma:\gamma < \kappa\}$ be a 
sequence of subsets of $2^\theta$ such that
$\gamma_1 \ne \gamma_2 \Rightarrow |X_{\gamma_1} \backslash X_{\gamma_2}| =
2^\theta$; let $\langle Y_j:j < \mu \rangle$ be a sequence of subsets of
$\kappa$ such that $j_1 \ne j_2 \Rightarrow |Y_{j_1} \backslash Y_{j_2}| = 
\kappa$, let
$g_\gamma$ be a one to one mapping from $\theta$ into $\dbcu_{n < \omega}
a^\gamma_n$ mapping $\theta_n$ onto $a^\gamma_n$, and lastly 
let $t^\gamma_{\ell,\alpha} = 
g''_\gamma(t_{\ell,\alpha})$ and $t^\gamma_{\ell,\alpha,\beta} =
g''_\gamma(t^\gamma_{\ell,\alpha,\beta})$.  Let $t^\gamma_{3,\alpha,\beta} =
\{g_\gamma(\varepsilon):\varepsilon \in t_{2,\alpha,\beta}$ and
$|t_{2,\alpha,\beta} \cap \varepsilon|$ is even$\}$.
Let $t^\gamma_{\zeta,\alpha,\beta} =: \{\zeta \in
t^\gamma_{2,\alpha,\beta}:|t^\gamma_{2,\alpha,\beta} \cap \zeta|$ is even$\}$.
\sn
For $j < \mu$, let ${\Cal A}_j$ be the following family of subsets of $\kappa$

$$
t^\gamma_{0,\alpha},t^\gamma_{1,\alpha} \text{ when } \gamma < \kappa,
\alpha < 2^\theta
$$

$$
t^\gamma_{2,\alpha,1 +\beta}  \text{ when } \beta \notin X_\gamma,\alpha <
2^\theta
$$

$$
t^\gamma_{3,\alpha,1 + \beta} \text{ when } \beta \in X_\gamma,\alpha <
2^\theta 
$$

$$
t^\gamma_{2,\alpha,0} \text{ when } \gamma \notin Y_j \text{ and }
t^\gamma_{3,\alpha,0} \text{ when } \gamma \in Y_j.
$$
\mn
Clearly
\mr
\item "{$\bigodot_1$}"  $t' \ne t'' \in {\Cal A}_j \Rightarrow |t' \cap t''|
< \aleph_0 = |t'|$.
\ermn
Let ${\Cal A}^+_j$ be a maximal almost disjoint family of countable subsets
of $\mu$ extending ${\Cal A}_j$.  Let $I_j$ be the Boolean ring of subsets
of $\kappa$ generated by ${\Cal A}^+_j \cup \{\{\varepsilon\}:\varepsilon <
\kappa\}$ and $B_j$ be the Boolean algebra of subsets of $\kappa$ generated by
$I_j$.  Now
\mr
\item "{$\bigodot_2$}"  if $i_0,i_1 < \mu$ and $b_0,b_1 \in [\kappa]^\theta$ 
and $h$ is a one to one mapping from $b_0$ onto $b_1$ such that $\alpha \in
\text{ Dom}(b) \Rightarrow h(\alpha) \ne \alpha$, \ub{then}
for some $t^0 \in {\Cal A}^+_{i_0},t^1 \in {\Cal A}^+_{i_1}$ we have:
$t^0 \subseteq^* b_0,t^1 \subseteq^* b_1$ and $h$ maps $t^0$ into a 
co-infinite subset of $t^1$ \nl
[why?  for some $\gamma_0 < \kappa$ the set $b_0 \cap \dbcu_{n < \omega}
a^{\gamma_0}_n$ have cardinality $\theta$, so without loss of generality
$b_0 \subseteq \dbcu_{n < \omega} a^{\gamma_0}_n$ and similarly for some
$\gamma_1 < \kappa$ without loss of generality
$b_1 \subseteq \dbcu_{n < \omega} a^{\gamma_1}_n$.  For $\ell =0,1$ let
$b^-_\ell \in [\theta]^\theta$ be such that $g_{\gamma^\ell}$ maps $b^-_\ell$
onto $b_\ell$.  Now \wilog \, $b^-_0 \cap b^-_1 = \emptyset$
(recall we have to preserve "$h$ is from $b_0$ onto $b_1$", too!).  
If $b^-_0 \cap b^-_1 = \emptyset$ 
then by clause $(v)$ of the fact some
$t^{\gamma^0}_{0,\alpha} \in {\Cal A}_{i_0} \subseteq {\Cal
A}^+_{i^0}$ and $t^{\gamma^1}_{0,\alpha_1} \in {\Cal A}_{i_1}
\subseteq {\Cal A}^+_{i_1}$ will be as
required in clause $(\alpha)$. 
So assume $b^-_0 = b^-_1$, let $b^*_0 = 
\{\alpha \in b^-_0:h \circ g_{\gamma_0}(\alpha) \ne g_{\gamma_1}
(\alpha)\}$.  If $b^*_0$ has cardinality $\theta$, we get the 
desired conclusion (in clause $(\alpha)$), so
assume $|b^*_0| < \theta$ hence \wilog \, $b^*_0 = \emptyset$.  
Also if $\gamma_0 \ne \gamma_1$ then
$|X_{\gamma_0} \backslash X_{\gamma_1}| = \kappa$ hence we can find a
non zero ordinal $\beta \in X_{\gamma_0} \backslash X_{\gamma_1}$ 
and we can find an ordinal $\alpha < 2^\theta$ such that $(\forall
\beta' < 2^\theta)[t^\gamma_{2,\alpha,\beta} \subseteq b^-_0]$
hence we can use
$t^\gamma_{3,\alpha,\beta},t^\gamma_{2,\alpha,\beta}$.  So we have to 
assume $\gamma_0 = \gamma_1$ but then $g_{\gamma_0} = g_{\gamma_1}$ so
$h \restriction (b_0 \backslash b^*_0)$ is the identity, a contradiction.]
\sn
\item "{$\bigodot_3$}"  if $i_0 \ne i_1$ and \footnote{by a little more care in
indexing, $Z \in [\mu]^{< \mu}$ is O.K. and we can choose $\gamma$
such that \nl
$\dbcu_n a_{\gamma,n} \le \kappa \backslash Z \backslash Z_0$}
$Z \in [\kappa]^{< \theta}$ and $h$ is a one to one function from
$\kappa \backslash Z$ onto $\kappa \backslash Z$ 
\ub{then} for some $t^0 \in {\Cal A}^+_{i_0},t^0 \subseteq^*
\text{ Dom}(h)$ and $t^1 \in {\Cal A}^+_{i_1}$ we have: $h''(t^0) \subseteq^* t^1,
t^1 \backslash h''(t^0)$ is infinite. \nl
[Why?  Let $Z_1 = \{\alpha \in \text{ Dom}(h):h(\alpha) \ne \alpha\}$, so by
$\odot_2$ we know $|Z_1| < \theta$.  We know that $Y_{i_0} \backslash Y_{i_1}$
has cardinality $\mu$, hence for some $\gamma \in Y_{i_0} \backslash Y_{i_1}$
we have set$(\bar a_\gamma) \cap [Z \cup Z_1]^{\aleph_0} = \emptyset$.  So
$t^\gamma_{3,\alpha,0} \in {\Cal A}_{i_0} \subseteq {\Cal A}^+_j$ and
$t^\gamma_{2,\alpha,0} \in {\Cal A}_{i_1} \subseteq {\Cal A}^+_{i_1}$, so
$t^\gamma_{3,\alpha,0}$ is a co-infinite subset of $t^\gamma_{2,\alpha,0},
t^\gamma_{2,\alpha,0} \subseteq^* \kappa \backslash Z \backslash Z_0$ and $h$
maps $t^\gamma_{3,\alpha,0} \backslash Z \backslash Z_0$ to itself a
co-infinite subset of $t^\gamma_{2,\alpha,0}$.]
\ermn
Clearly $\langle B_j:j < \mu \rangle$ is as required so we are done. \nl
2) Similar proof.  \hfill$\square_{\scite{ge.7}}$
\bigskip

\par \noindent \llap{---$\!\!>$} MARTIN WARNS: Label ge.8 on next line is also used somewhere else (Perhaps should have used scite instead of stag?\par
\demo{\stag{ge.8} Conclusion}  1) Under the assumption of \scite{ge.7},
let $\lambda^* = \text{ Ded}^+(\mu) = \text{ Min}\{\lambda$: there is
no tree with $\le \mu$ nodes and $\ge \lambda$ branches (equivalently,
a linear order of cardinality $\lambda$ and density $\le \mu\}$.
\ub{Then} fora ny $\lambda \in [\mu,\lambda^*)$ there is a superatomic
Boolean Algebra of cardinality $\lambda,\mu$ atoms with no
automorphism moving $\ge \theta$ atoms. \nl
2) Assume: $\theta$ is uncountable strong limit of cofinality
$\aleph_0$, pp$_{J^{\text{bd}}_\omega}(\theta) = 2^\theta$ (see
\cite[Ch.IX,\S5]{Sh:g} why this is reasonable) and $\kappa =
(2^\theta)^{+n},
\mu = 2^\kappa$ and $\mu < \lambda < \text{Ded}^+(\mu)$, 
e.g. $\lambda = 2^\chi$ for $\chi = \text{
Min}\{\chi:2^\chi > \mu\}$.  \ub{Then} there is a superatomic Boolean
Algebra of cardinality $\lambda$ and $\mu$ atoms, with no automorphism
moving $\ge \theta$ atoms. 
\enddemo
\bigskip

\demo{Proof}  Combine \scite{ge.7}, \scite{ge.6} and \scite{gd.2}.
\enddemo
\bigskip

\remark{Remark}  1)  So clearly in 
many models of ZFC we get that the bound is \scite{gd.1} cannot
be improved. \nl
2) The question is whether inductively we can get for many $\theta$'s the 
parallel of \scite{ge.7}. \nl
3) We can in \scitet{ge.8} replace $\theta$ by $\aleph_0$ (recall for
\scitet{ge.8}(2)) that we can replace the use of \scite{ge.6} by the known:
there is ${\Cal A} \subseteq [{}^\omega 2]^{\aleph_0}$ which is MAD,
every $A \subseteq {}^\omega 2$ dense in itself contains a member of
${\Cal A}$, each $A \in {\Cal A}$ has exactly one accumulation point.
\endremark
\newpage

\head {\S4} \endhead  \resetall \sectno=4
\bigskip

In the bound $\beth_4(\sigma)$, the last exponentiation was really
$sa(\mu)$ where
\definition{\stag{gf.1} Definition}  1) $sa^+(\mu) = 
\sup\{|B|^+:B \text{ is a superatomic Boolean Algebra with } \mu
\text{ atoms}\}$. \nl
2) $sa(\mu) = \sup\{|B|:B \text{ is a superatomic Boolean Algebra with } 
\mu \text{ atoms}\}$. \nl
3) $sa^+(\mu,\theta) = \sup\{|B|^+:B \text{ is a superatomic Boolean
Subalgebra of } {\Cal P}(\mu)$ extending $\{a \subseteq \mu:$ a finite
or cofinite such that $a \in B \Rightarrow (a) < \theta \vee |\mu
\backslash a| < \theta\}$. \nl
4) $sa(\mu,\theta) = \sup\{|B|:B \text{ is as in (3)}\}$. \nl
5) $sa^*(\theta) = \text{ Min}\{\lambda:\text{cf}(\lambda) \ge \theta$
and if $\mu < \lambda$ then $sa^+(\mu,\theta) \le \lambda\}$.
\enddefinition
\bn
That is, by the proof of Theorem \scite{gd.1}
\proclaim{\stag{gf.2} Claim}  If $B$ is a superatomic Boolean Algebra with no
automorphism moving $\ge \theta$ atoms, $\theta = \text{ cf}(\theta) >
\aleph_0$ then $|B| < sa^+(\beth_2(2^{< \theta})$, moreover $|B| < 
sa^+(\beth_2(sa^*(\theta))$.
\endproclaim
\bn
\ub{\stag{gf.3} Discussion}:  Now 
consistently $sa(\aleph_1) < 2^{\aleph_1}$, as \cite[8.1]{Sh:620} show
the consistency of a considerably stronger statement.  It proves that e.g.
if we start with $\bold V \models$ GCH and $\Bbb P$ is adding
$\aleph_{\omega_1}$ Cohen reals \ub{then} in $\bold V^{\Bbb P},
(2^{\aleph_0} = \aleph_{\omega_1} < 2^{\aleph_1} = \aleph_{\omega_1
+1}$ and) among any $\aleph_{\omega_1+1}$ members of ${\Cal P}(\omega_1)$
there are $\aleph_{\omega_1+1}$ which form an independent family, i.e. any
finite nontrivial Boolean combination of them is nonempty, in other words
``${\Cal P}(\omega_1)$ has $\aleph_{\omega_1+1}$-free precaliber in 
Monk's question
definition".  (Not surprising this is the same model for ``no tree with
$\aleph_1$ nodes has $2^{\aleph_1}$ branches" in \cite{B1}). \nl
So the bound $\beth_4(\theta)$ is not always the right ones.
\bigskip

\proclaim{\stag{gf.4} Claim}  Assume
\mr
\item "{$(a)$}"  $\Upsilon = \Upsilon^{< \Upsilon} < \mu = 
{\text{\rm cf\/}}(\mu) < \chi$
\sn
\item "{$(b)$}"  {\rm cf}$(\chi) = \mu$ and $(\forall \alpha < \chi)
(|\alpha|^\mu < \chi)$ and $(\forall \alpha < \mu)(|\alpha|^{< \Upsilon} <
\mu)$
\sn
\item "{$(c)$}"  $\Bbb Q$ is a forcing notion of cardinality $< \chi$ such that
in $\bold V^{\Bbb Q}:\mu$ is a regular cardinal $(\forall a \in [\chi]^{< \mu})(\exists b)
[a \subseteq b \in ([\chi]^{< \mu})^{\bold V}]$
\sn
\item "{$(d)$}"  $\Bbb P = \{f:f$ a partial function from $\chi$ to $\{0,1\}$ of
cardinality $< \Upsilon\}$ order by inclusion (that is, adding a $\chi \,
\sigma$-Cohen).
\ermn
\ub{Then} in $\bold V^{{\Bbb Q} \times {\Bbb P}}$ we 
have: $(2^\sigma = 2^{< \mu} = \chi,2^\mu =
\chi^\mu = (\chi^\mu)^{\bold V}$ and) $sa(\mu) = \chi < 2^\mu$, moreover the Boolean
Algebra ${\Cal P}(\mu)$ has $\chi^+$-free precaliber.
\endproclaim
\bigskip

\demo{Proof}  Work in $\bold V^{\Bbb Q}$, like \cite[8.1]{Sh:620}, not
using ``$\Bbb P$ is $\sigma$-complete" 
which may fail in $\bold V^{\Bbb Q}$. \hfill$\square_{\scite{gf.4}}$
\enddemo
\bn
On the other hand
\proclaim{\stag{gf.5} Claim}  Assume $\bar \lambda = \langle \lambda_n:n <
\omega \rangle$ satisfies $\lambda_{n+1} = \text{ Min}\{\lambda:2^\lambda >
2^{\lambda_n}\}$.  \ub{Then} for infinitely many $n$'s for some $\mu_n \in [
\lambda_n,\lambda_{n+1})$ we have $sa(\mu_n) = 2^{\mu_n} = 2^{\lambda_n}$
(in fact $sa^+(\mu_n) = (2^{\mu_n})^+ = (2^{\lambda_n})^+$ except possibly 
when cf$(2^{\lambda_n}) \le 2^{\lambda_{n-1}}$.
\endproclaim
\bigskip

\demo{Proof}  By \cite[3.4]{Sh:430} we have for infinitely many $n$'s
$\mu_n \in [\lambda_n,\lambda_{n+1})$ and for every regular $\chi \le
2^{\lambda_n} = 2^{\mu_n}$, a tree with $\le \mu$ nodes, $\lambda_n$ levels
and $\ge \chi$ \, $\lambda_n$-branches. \hfill$\square_{\scite{gf.5}}$
\enddemo
\bn
\ub{\stag{gf.6} Conclusion}:  1) Assume $\theta$ is strong limit, $\theta
> \text{ cf}(\theta) = \aleph_0$ and Pr$(2^{2^\theta},\theta)$ and $\lambda
< sa^+(\beth_3(\theta))$.  \ub{Then}
\mr
\item "{$(*)_{\theta,\lambda}$}"  there is a superatomic Boolean Algebra
without any automorphism moving $\ge \theta$ atoms such that $B$ has
cardinality $\lambda$ (and has $\lambda$ atoms \footnote{we can allow less
atoms and less elements}). 
\ermn
2) Assume Pr$(\beth_2,\aleph_1)$ and
$\lambda < sa^+(\beth_3)$.  Then $(*)_{\theta,\lambda}$ holds.
\bigskip

\demo{Proof}  1) Use \scite{ge.7} and \scite{gd.2}. \nl
2) Similar only replace \scite{ge.7} by a parallel claim. \nl
${{}}$  \hfill${\scite{gf.6}}$
\enddemo

\noindent
\centerline {$* \qquad * \qquad *$}
\mn
\par \noindent \llap{---$\!\!>$} MARTIN WARNS: Label ge.8 on next line is also used somewhere else (Perhaps should have used scite instead of stag?\par
\ub{\stag{ge.8} Discussion}:  [here?] Suppose 
$\theta = \tau^+$ and there is a tree
${\Cal T}$ with $\tau$ nodes and $\varrho > \theta$ branches.  We can build a
superatomic Boolean ring $B$ with $|{\Cal T}| \le \tau$ atoms and $\varrho$
elements, let $Y$ be a natural set of representatives for
the set of higher level atoms (i.e. not atoms).  By the
ulftrafilter for $y \in Y$ and say $B \subseteq {\Cal P}({\Cal T}),
\dsize \bigwedge_{t \in {\Cal T}} [\{t\} \in B]$ and $A$ a set of $\varrho$
elements $\{d_y:y \in Y_0 \subseteq Y\} \subseteq 
[A]^{\aleph_0}$ as in \S2, and use the
Boolean ring $B^\oplus$ generated by $\{\{t\}:t \in {\Cal T}\} \cup
\{y \cup d_y:y \in Y\}$.

We would like: every automorphism of $B^\oplus$ moves $\le \tau$ atoms.
If we succeed, we can continue to immitate \scite{ge.7} with the present
$\theta$.

So it is natural to consider:
\mr
\item "{$(*)_{\tau,\varrho_1,\varrho_2,\varrho_3}$}"  there is a 
tree with $\tau$ nodes and
a set of $\varrho_3$ branches such that an automorphism of the tree with moves
$> \tau$ branches move at least $\varrho_2$ branches.
\ermn
Restricting ourselves to $\varrho_0$-branches $\varrho_0 = 
\text{ cf}(\varrho_0) \le \varrho_2,2^{\varrho_0} = 2^\tau$, 
we can make the superatomic Boolean Algebras quite rigid so we need
\mr
\item "{$(*)_{\tau,\varrho_1,\varrho_2}$}"  there is a tree ${\Cal T}$ 
with $\tau$ nodes such that any subtree has $\le \varrho_1$ and 
$\ge \varrho_2$ branches.
\ermn
\bigskip

\newpage
    
REFERENCES.  
\bibliographystyle{lit-plain}
\bibliography{lista,listb,listx,listf,liste}

\shlhetal
\enddocument

\bye